\documentclass{article}
\usepackage[T2A]{fontenc}
\usepackage[cp1251]{inputenc}
\usepackage[english,russian]{babel}
\usepackage[tbtags]{amsmath}
\usepackage{amsfonts,amssymb,mathrsfs,amscd}
\numberwithin{equation}{section}
\newtheorem{remark}{Remark}[section]
\newtheorem{lemma}{Lemma}[section]
\newtheorem{theorem}{Theorem}[section]
\newtheorem{definition}{Definition}[section]

\newtheorem{conclusion}{Conclusion}
\newtheorem{problem}{Problem}[section]

\renewcommand{\Re}{\mathop{\rm Re}}

\newcommand{\tr}{\mathop{\rm tr}}

\newcommand{\col}{\mathop{\rm col}}

\begin{document}
\begin{Large}
\thispagestyle{empty}
\begin{center}
{\bf Inverse scattering problem for the fifth-order Schr${\rm\bf\ddot{o}}$dinger operator on the whole axis\\
\vspace{5mm}
V. A. Zolotarev}\\

ORCID ID 0000-0002-2165-8552

B. Verkin Institute for Low Temperature Physics and Engineering
of the National Academy of Sciences of Ukraine\\
47 Nauky Ave., Kharkiv, 61103, Ukraine

\end{center}
\vspace{5mm}

{\small {\bf Abstract.} A method of solution for the inverse scattering problem for the linear fifth order Schr${\rm\ddot{o}}$dinger operator is found. Two closed systems of linear singular integral equations are obtained. From solutions to these equations, potential on the left and right half-axes is restored unambiguously.}
\vspace{5mm}


\vspace{5mm}

\begin{center}
{\bf Introduction}
\end{center}
\vspace{5mm}

Theory of inverse problems is one of the most interesting and intriguing fields of mathematical physics. Inverse scattering problems originating in quantum mechanics take a special place in it. Inverse scattering problems for the classical Scr${\rm\ddot{o}}$dinger operator were solved by V. A. Marchenko \cite{1} and M. G. Krein \cite{2}. The methods of works \cite{1,2} were further developed in the studies of V. A. Marchenko \cite{3}, L. D. Faddeev \cite{4}, Deift -- Trubovitz \cite{5}. The method of inverse scattering problem found its effective application in the problems of integration of nonlinear equations (Korteveg -- de Vries, sin-Gordon, etc.). Solution of the inverse scattering problem for differential operators of the higher ($>2$) orders turns out to be not so simple and it is caused, first of all, by the absence of transformation operators \cite{3,6} (when $n>2$). And in the cases when such operators exist, they give a rather narrow class of potentials consisting of analytic functions. Inverse scattering problems for third-order operators are solved in \cite{7,8}.

Many equations of mathematical physics (e.g., Kaup -- Kupershmidt, Sawada -- Kotega, Candrey -- Dodd -- Gibbon, etc.) are of fifth order. Moreover, operators describing neural networks are of fifth order \cite{10,12}. One has to mention that in the works \cite{9} -- \cite{16} there emerge fifth-order equations. Therefore, the inverse scattering problem for the fifth-order Schr${\rm\ddot{o}}$dinger operator emerged naturally.

This work studies the inverse scattering problem for the fifth-order Schr${\rm\ddot{o}}$dinger equation
\begin{equation}
\left(\frac1iD\right)^5y(\lambda,x)+q(x)y(\lambda,x)=\lambda^5y(\lambda,x)\quad(x\in\mathbb{R})\label{eq0.1}
\end{equation}
where ${\displaystyle D=\frac d{dx}}$, $\lambda\in\mathbb{C}$, and $q(x)$ is a real function ({\bf potential}) such that
\begin{equation}
\int\limits_{\mathbb{R}}e^{2a|x|}|q(x)|^2dx<\infty\quad(a>0).\label{eq0.2}
\end{equation}
Operator $L_q{\displaystyle=\left(\frac1iD\right)^5+q(x)}$ with domain $\mathfrak{D}(L_q)=W_2^5(\mathbb{R})$ is self-adjoint.

The paper consists of 4 sections. Section 1 describes the fundamental system of solutions to the unperturbed equation \eqref{eq0.1} (where $q(x)\equiv0$). Here, following \cite{6,7}, functions $\{s_k(i\lambda x)\}_0^4$ are defined (analogues of classical sines and cosines). System $\{s_k(i\lambda x)\}_0^2$ for the third-order equation \cite{7,8} turned out to be the instrument that allowed inverse scattering problems for third-order operators to be solved. This section lists the basic properties of the system of functions $\{s_k(i\lambda x)\}_0^4$ which deserves special attention.

Second section constructs Jost solutions to equation \eqref{eq0.1}. Besides the proof of the existence of Jost solutions, which is due to \eqref{eq0.2}, here also analytic properties with regard to the spectral parameter $\lambda$ of such Jost solutions are studied. In particular, analyticity of Jost solutions in proper sectors is established. Explicit formulas allowing to recover potential $q(x)$ on the left and right half-axes from Jost solutions are obtained.

The main section of the work is Section 3 which studies the scattering problem for waves incident from ``$+\infty$''. Wronskian of a part of Jost solutions and its analytic properties is an important instrument of the study. Using properties of such partial Wronskian, the basic properties of the transition matrix are obtained. These properties characterize unitarity of the scattering problem. Further, using the method of Riemann boundary value problems, the closed linear system of singular integral equations of inverse problem is derived. It is shown how one can recover potential $q(x)$ on the right half-axis, $x\in\mathbb{R}_+$, from the solution to this system.

The last Section 4 studies the scattering problem for the waves incident from ``$-\infty$''. As in Section 3, here the closed system of linear singular integral equations of the inverse problem is obtained. From the solution to this system, potential $q(x)$ is restored unambiguously on the left half-axis $x\in\mathbb{R}_-$.

So, we obtained two systems of linear singular integral equations, the solutions to which allow to completely recover potential $q(x)$.

\section{Preliminary information}\label{s1}

{\bf 1.1} Solution to the fifth-order equation
\begin{equation}
\left(\frac1iD\right)^5y(x)=\lambda^5y(x)\quad(x\in\mathbb{R},\,\lambda\in\mathbb{C},\,D=d/dx)\label{eq1.1}
\end{equation}
has five linear independent solutions $\{e^{i\lambda\zeta_kx}\}_0^4$ where $\{\zeta_k\}_0^4$ are the roots of equation $z^5=1$,
\begin{equation}
\zeta_k=\exp\left\{\frac{i2\pi}5\right\}_{k=0}^4,\label{eq1.2}
\end{equation}
besides, it is obvious that
\begin{equation}
\zeta_0=1;\quad\overline{\zeta_1}=\zeta_4;\quad\overline{\zeta_2}=\zeta_3;\quad\zeta_k=\zeta_1^k\quad(0\leq k\leq4).\label{eq1.3}
\end{equation}
Another system of fundamental solutions $\{s_p(i\lambda x)\}$ to equation \eqref{eq1.1} plays an important role. This system for the third-order equation is used effectively in \cite{7,8} where
\begin{equation}
s_p(z)\stackrel{\rm def}{=}\frac15\sum\limits_k\frac1{\zeta_k^p}e^{\zeta_kz}\quad(0\leq p\leq4).\label{eq1.4}
\end{equation}

\begin{remark}
Definitions \eqref{eq1.4} and \eqref{eq1.3} imply

{\rm(i)} {\bf reality}, $\overline{s_p(z)}=s_p(\overline{z})$ ($0\leq p\leq 4$);

{\rm(ii)} {\bf $p$-evenness}, $s_p(\zeta,z)=\zeta_1^ps_p(z)$ ($0\leq p\leq4$);

{\rm(iii)} ${\displaystyle\frac d{ds}s_p(z)=s_{p-1}(z)}$ ($1\leq p\leq4$); ${\displaystyle\frac d{dz}s_0(z)=s_4(z);}$

{\rm(iv)} functions $\{s_p(z)\}_0^4$ are solutions to equation $y^{(5)}(z)=y(z)$ and satisfy the initial data
\begin{equation}
\begin{array}{cccc}
s_0(0)=1,&s'(0)=0,&...&s^{(4)}(0)=0;\\
s_1(0)=1,&s'_1(0)=1,&...&s_1^{(4)}(0)=0;\\
...&...&...&...\\
s_4(0)=0,&s'_4(0)=0,&...&s_4^{(4)}(0)=1.
\end{array}\label{eq1.5}
\end{equation}
\end{remark}

By the functions $\{s_p(z)\}_0^4$ \eqref{eq1.4}, we construct the matrix-function
\begin{equation}
\mathbb{S}(z)\stackrel{\rm def}{=}\left[
\begin{array}{cccc}
s_0(z)&s_1(z)&...&s_4(z)\\
s_4(z)&s_0(z)&...&s_3(z)\\
...&...&...&...\\
s_1(z)&s_2(z)&...&s_0(z)
\end{array}\right].\label{eq1.6}
\end{equation}
Equations (iii) \eqref{eq1.5} and (iv) \eqref{eq1.5} imply that $\mathbb{S}(z)$ is the solution to Cauchy problem
\begin{equation}
\frac d{dz}\mathbb{S}(z)=\tau\mathbb{S}(z);\quad\mathbb{S}(0)=I.\label{eq1.7}
\end{equation}
where
\begin{equation}
\tau\stackrel{\rm def}{=}\left[
\begin{array}{ccccc}
0&1&0&...&0\\
0&0&1&...&0\\
...&...&...&...&...\\
1&0&0&...&0
\end{array}\right]\quad(\tau^5=I).\label{eq1.8}
\end{equation}
Note that the correspondence $\zeta_1\rightarrow\tau$ is the matrix representation of a cyclic group $\{\zeta_1^k\}_0^4$ \eqref{eq1.2} of identity roots.

\begin{lemma}\label{l1.1}
For the functions $\{s_p(z)\}_0^4$ \eqref{eq1.4}, the following statements are true:

{\rm(i)} the {\bf Euler formula}
$$e^{\zeta_kz}=s_0(z)+\zeta_ks_1(z)+...+\zeta_k^4s_4(z)\quad(0\leq k\leq4);$$

{\rm(ii)} the {\bf main identity}
$$\det\mathbb{S}(z)=1$$
($\mathbb{S}(z)$ is given by \eqref{eq1.6});

{\rm(iii)} {\bf addition formulas}
$$\mathbb{S}(z+\zeta_1^kw)-\mathbb{S}(z)(b_\zeta^*)^k\mathbb{S}(w)b_\zeta^k\quad(\mathbb{S}(z+w)=\mathbb{S}(z)\cdot \mathbb{S}(w),\,k=0);$$
where
$$b_\zeta\stackrel{\rm def}{=}\left[
\begin{array}{cccc}
1&0&...&0\\
0&\zeta_1&...&0\\
...&...&...&...\\
0&0&...&\zeta_4
\end{array}\right]\quad(b_\zeta^*b_\zeta=I).$$

{\rm(iv)} $p$-evenness of $\mathbb{S}(z)$ \eqref{eq1.6},
$$\mathbb{S}(z\zeta_1)=b_\zeta^*\mathbb{S}(z)b_\zeta;$$

{\rm(v)}
\begin{equation}
\widehat{s}_{0,0}=s_{0,0}(-z);\quad\widehat{s}_{0,1}(z)=s_{0,1}(-z),...,\widehat{s}_{0,4}(z)=s_{0,4}(-z);\label{eq1.9}
\end{equation}
where $\widehat{s}_{k,s}(z)$ are algebraic compliments of the elements $s_{k,s}(z)$ of the matrix $\mathbb{S}(z)$ \eqref{eq1.6};

{\rm(vi)} {\bf Taylor formulas}
$$s_0(z)=I+\frac{z^5}{5!}+\frac{z^{10}}{10!}+...;$$
$$s_1(z)=\frac z{1!}+\frac{z^6}{6!}+\frac{z^{11}}{11!}+...;$$
$$s_2(z)=\frac{z^2}{2!}+\frac{z^7}{7!}+\frac{z^{12}}{12!}+...;$$
$$s_3(z)=\frac{z^3}{3!}+\frac{z^8}{8!}+\frac{z^{13}}{13!}+...;$$
$$s_4(z)=\frac{z^4}{4!}+\frac{z^9}{9!}+\frac{z^{14}}{14!}+... .$$
\end{lemma}

P r o o f. Equality (i) \eqref{eq1.9} follows from the definition \eqref{eq1.4} of the functions $\{s_p(z)\}_0^4$. In view of the Cauchy problem \eqref{eq1.7}, we have
\begin{equation}
\mathbb{S}(z)=\exp\{z\tau\}.\label{eq1.10}
\end{equation}
Further using the well-known equality $\det(\exp B)=\exp(\tr B)$ and the fact that $\tr\tau=0$, we obtain that $\det\mathbb{S}(z)=1$. Note that \eqref{eq1.10} implies that
$$\mathbb{S}(z+w)=\mathbb{S}(z)\mathbb{S}(w).$$
Equality (iv) \eqref{eq1.9} follows from $p$-evenness (ii) \eqref{eq1.5} of the functions $s_p(z)$. Note that relations (v) \eqref{eq1.9} are a corollary of equality $\mathbb{S}^{-1}=\mathbb{S}(-z)$. Taylor formulas (vi) \eqref{eq1.9} follow from \eqref{eq1.10} and the fact that $\tau^5=I$. $\blacksquare$

\begin{remark}\label{r1.2}
The main identity {\rm(ii)} \eqref{eq1.9} in the two-dimensional case for
$$\mathbb{S}(z)=\left[
\begin{array}{ccc}
\cosh z&\sinh z\\
\sinh z&\cosh z
\end{array}\right]$$
is $\cosh^2z-\sinh^2z=1$. Transformation $\widetilde{h}=s(z)h$ ($h=\col(x,y)$, $\widetilde{h}=\col(\widetilde{x},\widetilde{y})$), i.e., $\widetilde{x}=x\cosh z+y\sinh z$, $\widetilde{y}=x\sinh z+y\cosh z$, in the matrix form becomes
$$\left[
\begin{array}{ccc}
\widetilde{x}&\widetilde{y}\\
\widetilde{y}&\widetilde{x}
\end{array}\right]=\mathbb{S}(z)\left[
\begin{array}{ccc}
x&y\\
y&x
\end{array}\right].$$
If $\det s(z)=1$, then transformation $\widetilde{h}=\mathbb{S}(z)h$ preserves quadratic forms $\widetilde{x}^2-\widetilde{y}^2=x^2-y^2$.

In the three-dimensional situation
\begin{equation}
\mathbb{S}(z)=\left[
\begin{array}{ccc}
s_0(z)&s_1(z)&s_2(z)\\
s_2(z)&s_0(z)&s_1(z)\\
s_1(z)&s_2(z)&s_0(z)
\end{array}\right],\label{eq1.12}
\end{equation}

The main identity is
\begin{equation}
s_0^3(z)+s_1^3(z)+s_2^3(z)-3s_0(z)s_1(z)s_2(z)=1\quad(\det\mathbb{S}(z)=1)\label{eq1.12}
\end{equation}
Similarly, transformation $\widetilde{h}=\mathbb{S}(z)h$ ($h=\col(x,y,z)$, $\widetilde{h}=\col(\widetilde{x},\widetilde{y},\widetilde{z})$) given by $\mathbb{S}(z)$ \eqref{eq1.12}
$$\widetilde{x}=xs_0(z)+ys_1(z)+zs_2(z);\quad\widetilde{y}=xs_2(z)+ys_0(z)+zs_1(z);$$
$$\widetilde{z}=xs_1(z)+ys_2(z)+zs_0(z)$$
which in the matrix form becomes
$$\left[
\begin{array}{ccc}
\widetilde{x}&\widetilde{z}&\widetilde{y}\\
\widetilde{y}&\widetilde{x}&\widetilde{z}\\
\widetilde{z}&\widetilde{y}&\widetilde{x}
\end{array}\right]=\mathbb{S}(z)\left[
\begin{array}{ccc}
x&z&y\\
y&x&z\\
z&y&x
\end{array}\right],$$
and thus
$$\widetilde{x}^3+\widetilde{y}^3+\widetilde{z}^3-3\widetilde{x}\widetilde{y}\widetilde{z}=\det\mathbb{S}(z)[x^3+y^3+z^3-3xyz].$$
Therefore, identity $\det s(z)=1$ implies that transformation $h\rightarrow\widetilde{h}=\mathbb{S}(z)h$ preserves the cubic form $x^3+y^3+z^3-3xyz$ which is equal to the volume of parallelepiped formed by the vectors $h=\col(x,y,z)$, $\widehat{h}=\col(z,x,y)$, and $\widehat{\widehat{h}}=\col(y,z,x)$.

To obtain form of the main identity in the five-dimensional case, we note that $\mathbb{S}(z)$ \eqref{eq1.6} can be represented as
\begin{equation}
\mathbb{S}(z)=s_0(z)I+s_1(z)\tau+s_2(z)\tau^2+s_3(z)\tau^3+s_4(z)\tau^4\label{eq1.14}
\end{equation}
where $\tau$ is given by \eqref{eq1.8}. Equation \eqref{eq1.14} implies that the function $\det\mathbb{S}(z)=F(s_0(z),s_1(z),s_2(z),s_3(z),s_4(z))$ is invariant relative to the group of permutations $\{P^k\}$ where
\begin{equation}
P(s_0,s_1,s_2,s_3,s_4)=(s_4,s_0,s_1,s_2,s_3).\label{eq1.15}
\end{equation}
Hence, we obtain the main identity for $p=5$,
$$s_0^5(z)+s_1^5(z)+s_2^5(z)+s_3^5(z)+s_4^5(z)-5s_0(z)s_1(z)s_2(z)s_3(z)s_4(z)$$
\begin{equation}
+5\sum\limits_{k=0}^4p^k\{s_0(z)[s_2^2(z)s_3^2(z)+s_1^2(z)s_4^2(z)]\}-5\sum\limits_{k=0}^4p^k\{s_0^3(z)[s_2(z)s_3(z)\label{eq1.16}
\end{equation}
$$+s_1(z)s_4(z)]\}=1.$$
Besides, summation in \eqref{eq1.16} is over all the powers of the permutation $p$ \eqref{eq1.15}. As in the case of $p=3$, transformation $s(z)$ \eqref{eq1.6} preserves the form $F(x,y,z,u,v)$, where
$$F(x,y,z,u,v)=x^5+y^5+z^5+u^5+v^5-5xyzuv$$
$$+5\sum\limits_{k=0}^4p^k\{x[y^2v^2+z^2u^2]\}-5\sum\limits_{k=0}^4p^k\{x^3[yv+zu]\}.$$
\end{remark}
\vspace{5mm}

{\bf 1.2} By $L_{\zeta_k}$, we denote the straight lines in $\mathbb{C}$ with the direction vectors $\{\zeta_k\}_0^4$ \eqref{eq1.2},
\begin{equation}
L_{\zeta_k}\stackrel{\rm def}{=}\{x\zeta_k:x\in\mathbb{R}\}\quad(0\leq k\leq4).\label{eq1.17}
\end{equation}
And let $l_{\zeta_k}$ be rays in the direction of $\zeta_k$ exiting the origin and $\widehat{l}_{\zeta_k}$ be rays in the direction of $\zeta_k$ coming to the origin,
\begin{equation}
l_{\zeta_k}\stackrel{\rm def}{=}\{x\zeta_k:x\in\mathbb{R}_+\};\quad\widehat{l}_{\zeta_k}\stackrel{\rm def}{=}L_{\zeta_k}\backslash l_{\zeta_k}=\{x\zeta_k:x\in\mathbb{R}_-\}\quad(0\leq k\leq4).\label{eq1.18}
\end{equation}
It is obvious that
\begin{equation}
\zeta_1L_{\zeta_k}=L_{\zeta_{k+1}}\quad(0\leq k\leq3);\quad\zeta L_{\zeta_4}=L_{\zeta_0}.\label{eq1.19}
\end{equation}
Straight lines $\{L_{\zeta_k}\}_0^4$ \eqref{eq1.17} divide plane $\mathbb{C}$ into 10 sectors,
\begin{equation}
S_l\stackrel{\rm def}{=}\left\{z\in\mathbb{C}:\frac{2\pi}{10}l<\arg z<\frac{2\pi}{10}(l+1)\right\}\quad(0\leq l\leq9).\label{eq1.20}
\end{equation}

\begin{picture}(200,200)
\put(0,100){\vector(1,0){200}}
\put(100,0){\vector(0,1){200}}
\put(50,200){\vector(1,-2){100}}
\put(50,0){\vector(1,2){100}}
\put(200,67){\vector(-3,1){200}}
\put(200,133){\vector(-3,-1){200}}
\put(135,150){$S_1$}
\put(102,190){$S_2$}
\put(150,190){$L_{\zeta_1}$}
\put(154,105){$S_0$}
\put(190,105){$L_{\zeta_0}$}
\put(40,33){$S_6$}
\put(60,5){$S_7$}
\put(120,5){$L_{\zeta_4}$}
\put(20,105){$S_4$}
\put(20,85){$S_5$}
\put(10,60){$L_{\zeta_3}$}
\put(170,35){$S_8$}
\put(170,79){$S_9$}
\put(45,150){$S_3$}
\put(15,140){$L_{\zeta_2}$}
\qbezier(123,109)(100,130)(82,106)
\qbezier(100,80)(120,90)(120,107)
\qbezier(80,107)(80,80)(100,77)
\end{picture}

\hspace{20mm} Fig. 1

Straight line $L_{\zeta_0}$ divides sectors $S_0$ and $S_9$ and also sectors $S_4$ and $S_5$; straight line $L_{\zeta_1}$ divides sectors $S_1$ and $S_2$ and also $S_6$ and $S_7$. Similarly, straight line $L_{\zeta_2}$ divides the sectors $S_3$ and $S_4$ and also $S_0$ and $S_1$. And the straight line $L_{\zeta_3}$ divides the sectors $S_5$ and $S_6$ and also $S_0$ and $S_1$. Finally, $L_{\zeta_4}$ divides the sectors $S_7$ and $S_8$, and also the sectors $S_3$ and $S_2$.

Note also that $\widehat{l}_{\zeta_3}$ is the bisectrix of the angle between $l_{\zeta_0}$ and $l_{\zeta_1}$, and $\widehat{l}_{\zeta_4}$ is the bisectrix of the angle between $l_{\zeta_1}$ and $l_{\zeta_2}$. Similarly, $\widehat{l}_{\zeta_0}$ is the bisectrix of the angle between $l_{\zeta_2}$ and $l_{\zeta_3}$, and, finally, $l_{\zeta_1}$ is the bisectrix between $l_{\zeta_3}$ and $l_{\zeta_4}$.

\begin{problem}\label{p1.1}
In the works \cite{7,8}, it is established that zeros of the functions $s_0(z)$, $s_1(z)$, $s_2(z)$ ($p=3$) lie on the bisectrices of the corresponding sectors and are simple. Describe the set of zeros of the functions $\{s_p(z)\}_0^4$ ($p=5$), show where they are situated and characterize their multiplicity.
\end{problem}
\vspace{5mm}

{\bf 1.3} Consider the Cauchy problem
\begin{equation}
\begin{array}{ccc}
{\displaystyle\left(\frac1iD\right)^5y(x)=\lambda^5y(x)+f(x);}\\
y(0)=y_0,\,y'(0)=y_1,\,y''(0)=y_2,\,y'''(0)=y_3,\,y^{\rm iv}(0)=y_4
\end{array}\label{eq1.21}
\end{equation}
($\lambda\in\mathbb{C}$, $x\in\mathbb{R}$, $y_k\in\mathbb{C}$ ($0\leq k\leq4$)). Solution to the homogenous ($f(x)\equiv0$) problem, in view of (iv) \eqref{eq1.5}, is
\begin{equation}
y_0(\lambda,x)=\sum\limits_{k=0}^4y_k\frac{s_k(i\lambda x)}{(i\lambda)^k}.\label{eq1.22}
\end{equation}
Hence, by the method of variation of parameters, we find the general ($f(x)\not\equiv0$) solution to Cauchy problem \eqref{eq1.21}
\begin{equation}
y(\lambda,x)=\sum\limits_{k=0}^4y_k\frac{s_k(i\lambda x)}{(i\lambda)^k}+i\int\limits_0^x\frac{s_4(i\lambda(x-t))}{(i\lambda)^4}f(t)dt.\label{eq1.23}
\end{equation}

\section{Jost solutions}\label{s2}

{\bf 1.2} In the space $L^2(\mathbb{R})$, we consider the fifth-order linear operator
\begin{equation}
(L_qy)(x)=\left(\frac1iD\right)^5y(x)+q(x)y(x)\label{eq2.1}
\end{equation}
with domain $\mathfrak{D}(L_q)=W_2^5(\mathbb{R})$, assuming that potential $q(x)$ is a real function such that
\begin{equation}
\int e^{2a|x|}|q(x)|^2dx<\infty\quad(0<a<\infty).\label{eq2.2}
\end{equation}

\begin{remark}\label{r2.1}
Condition \eqref{eq2.2} implies that $e^{b|x|}|q(x)|\in L^1(\mathbb{R})\cap L^2(\mathbb{R})$, for all $b$ such that $0\leq b<a$.
\end{remark}

Consider equation
\begin{equation}
\left(\frac1iD\right)^5y(x)+q(x)y(x)=\lambda^5y(x)\label{eq2.3}
\end{equation}
which, in view of condition \eqref{eq2.2}, for $|x|\rightarrow\infty$, transforms into equation \eqref{eq1.1}; therefore, it is natural to define the Jost solutions $\{v_k(\lambda,x)\}_0^4$ and $\{u_k(\lambda,x)\}_0^4$ which are solutions to equation \eqref{eq2.3} and satisfy the boundary conditions
\begin{equation}
\begin{array}{ccc}
(a)\,v_k(\lambda,x)\rightarrow e^{i\lambda\zeta_kx}\quad(x\rightarrow\infty,\,0\leq k\leq4);\\
(b)\,u_k(\lambda,x)\rightarrow e^{i\lambda\zeta_kx}\quad(x\rightarrow-\infty,\,0\leq k\leq4).
\end{array}\label{eq2.4}
\end{equation}
Show that the integral equation
\begin{equation}
v_k(\lambda,x)=e^{i\lambda\zeta_kx}+i\int\limits_x^\infty\frac{s_4(i\lambda(x-t))}{(i\lambda)^4}q(t)v_k(\lambda,t)dt\quad(0\leq k\leq4)\label{eq2.5}
\end{equation}
is equivalent to the boundary value problem \eqref{eq2.3}, (a) \eqref{eq2.4} for $v_k(\lambda,x)$ ($0\leq
 k\leq4$). The fact that solution $v_k(\lambda,x)$ to problem (a) \eqref{eq2.4} one ought to seek in such form follows from \eqref{eq1.22}, \eqref{eq1.23}. In $L^2(\mathbb{R})$, we define the set of Volterra operators
$$(K_\lambda f)(x)=\int\limits_x^\infty K_1(\lambda,x,t)q(t)f(t)dt$$
where
\begin{equation}
K_1(\lambda,x,t)=\frac{s_4(i\lambda(x-t))}{(i\lambda)^4}.\label{eq2.6}
\end{equation}
In terms of $K_\lambda$, equation \eqref{eq2.5} becomes
$$(I-iK_\lambda)v_k(\lambda,x)=e^{i\lambda\zeta_kx},$$
and thus
\begin{equation}
v_k(\lambda,x)=\sum\limits_{n=0}^\infty i^nK_\lambda^ne^{i\lambda\zeta_kx}.\label{eq2.7}
\end{equation}
Operators $K_\lambda^n$ are also Volterra,
$$(K_\lambda^nf)(x)=\int\limits_x^\infty K_n(\lambda,x,t)q(t)f(t)dt,$$
and for their kernels $K_n(\lambda,x,t)$ the recurrent relations hold:
\begin{equation}
K_{n+1}(\lambda,x,t)=\int\limits_x^tK_n(\lambda,x,s)q(s)K_1(\lambda,s,t)ds.\label{eq2.8}
\end{equation}

\begin{lemma}\label{l2.1}
For the kernels $K_n(\lambda,x,t)$, the following estimates hold:
\begin{equation}
|K_n(\lambda,x,t)|\leq\frac{e^{|\lambda|(t-x)}}{|\lambda|^{4n}}\frac1{(n-1)!}\left\{\int\limits_x^t|q(s)|ds\right\}^{n-1}\quad(t>x>0).\label{eq2.9}
\end{equation}
\end{lemma}

P r o o f. For $K_1(\lambda,x,t)$ \eqref{eq2.6}, we obtain
$$|K_1(\lambda,x,t)|\leq\frac{e^{|\lambda|(t-x)}}{|\lambda|^4}\quad(t>x),$$
which coincides with \eqref{eq2.9}, for $n=1$. Using induction over $n$, equality \eqref{eq2.8} and estimates \eqref{eq2.9}, we have
$$|K_{n+1}(\lambda,x,t)|\leq\int\limits_x^t\frac{e^{|\lambda|(s-x)}}{|\lambda|^{4n}}\frac{\left\{\displaystyle\int\limits_x^s|q(\mu)|d\mu\right\}^{n-1}}{(n-1)!}\frac{|q(s)|e^{|\lambda|(t-s)}}{|\lambda|^4}ds$$
$$=\frac{e^{|\lambda|(t-x)}}{|\lambda|^{4(n+1)}}\frac1{n!}\left\{\int\limits_x^t|q(s)|ds\right\}^n.\blacksquare$$
Equation \eqref{eq2.7} implies that
\begin{equation}
v_k(\lambda,x)=e^{i\lambda\zeta_kx}+i\int\limits_x^\infty N(\lambda,x,t)q(t)e^{i\lambda\zeta_kt}dt\quad(0\leq k\leq4)\label{eq2.10}
\end{equation}
where
\begin{equation}
N(\lambda,x,t)\stackrel{\rm def}{=}\sum\limits_{n=1}^\infty i^{n-1}K_n(\lambda,x,t).\label{eq2.11}
\end{equation}
Series \eqref{eq2.11}, in view of \eqref{eq2.9}, is majorized by the convergent series,
\begin{equation}
|N(\lambda,x,t)|\leq\frac{e^{|\lambda|(t-x)}}{|\lambda|^4}\exp\left\{\frac1{|\lambda|^4}\int\limits_x^t|q(s)|ds\right\}.\label{eq2.12}
\end{equation}
Therefore, series \eqref{eq2.11} converges uniformly in $x$, which provides the needed differentiability with respect to $x$ of the function $v_k(\lambda,x)$.

Proceed to condition (a) \eqref{eq2.4}. Taking into account \eqref{eq2.10}, \eqref{eq2.12}, we obtain
$$|v_k(\lambda,x)-e^{i\lambda\zeta_kx}|\leq\int\limits_x^\infty\frac{e^{|\lambda|(t-x)}}{|\lambda|^4}\exp\left\{\frac1{|\lambda|^4}\int\limits_x^t|q(s)|ds\right\}\cdot|q(t)|e^{|\lambda|t}.$$
Since
$$\int\limits_x^t|q(s)|ds<\sigma_+,\quad\sigma_+\stackrel{\rm def}{=}\int\limits_0^\infty|q(s)|ds,$$
hence it follows that
$$|v_k(\lambda,x)-e^{i\lambda\zeta_kx}|\leq\frac{\exp\left\{\displaystyle\frac{\sigma_+}{|\lambda|^4}\right\}}{|\lambda|^4}\int\limits_x^\infty e^{|\lambda|(2t-x)}|q(t)dt.$$
Taking into account that $e^{|\lambda|(2t-x)}|q(t)|\leq e^{2|\lambda|t}|q(t)|\in L^1(\mathbb{R})\cap L^2(\mathbb{R})$ (Remark \ref{r2.1}), for $2|\lambda|<a$, we obtain that the last integral for $2|\lambda|<a$ vanishes when $x\rightarrow\infty$. This proves (a) \eqref{eq2.4}. Define the disk
\begin{equation}
\mathbb{D}_{a/2}=\{\lambda\in\mathbb{C}:|\lambda|<a/2\}.\label{eq2.13}
\end{equation}

\begin{theorem}\label{t2.1}
If \eqref{eq2.2} is true, then for all $\lambda\in\mathbb{D}_{a/2}$ \eqref{eq2.13}, Jost solutions $\{v_k(\lambda,x)\}_0^4$ of the boundary value problem \eqref{eq2.3}, (a) \eqref{eq2.4} exist and are given by \eqref{eq2.10}. Functions $\{v_k(\lambda,x)\}_0^4$ are holomorphic (in $\lambda$) in the disk $\mathbb{D}_{a/2}$.
\end{theorem}
\vspace{5mm}

{\bf 2.2} Analogously to \eqref{eq2.5}, equations
\begin{equation}
u_k(\lambda,x)=e^{i\lambda\zeta_k}-i\int\limits_{-\infty}^x\frac{s_k(i\lambda(x-t))}{(i\lambda)^4}q(t)u_k(\lambda,t)dt\quad(0\leq k\leq4)\label{eq2.14}
\end{equation}
correspond to Jost solutions $\{u_k(\lambda,x)\}_0^2$. In $L^2(\mathbb{R})$, we define the family of Volterra operators
$$(\widehat{K}_\lambda f)(x)=\int\limits_{-\infty}^xK_1(\lambda,x,t)q(t)f(t)dt$$
where $K_1(\lambda,x,t)$ is given by \eqref{eq2.6}, then equation \eqref{eq2.14} becomes
$$(I+i\widetilde{K}_\lambda)u_k(\lambda,x)=e^{i\lambda\zeta_kx}\quad(0\leq k\leq4),$$
and thus
\begin{equation}
u_k(\lambda,x)=\sum\limits_{n=0}^\infty(-i)^n\widehat{K}_\lambda^ne^{i\lambda\zeta_kx}.\label{eq2.15}
\end{equation}
Operators $\widehat{K}_\lambda^n$ are Volterra,
$$(\widehat{K}_\lambda^nf)(x)=\int\limits_{-\infty}^xK_n(\lambda,x,t)q(t)f(t)dt,$$
and for their kernels $\widehat{K}_n(\lambda,x,t)$, analogously to \eqref{eq2.9},
\begin{equation}
|\widehat{K}_n(\lambda,x,t)|\leq\frac{e^{|\lambda|(x-t)}}{|\lambda|^{4n}}\frac1{(n-1)!}\left\{\int\limits_t^x|q(s)|ds\right\}^{n-1}\quad(0>x>t).\label{eq2.16}
\end{equation}
Equation \eqref{eq2.15} yields that
\begin{equation}
u_k(\lambda,x)=e^{i\lambda\zeta_kx}+\int\limits_{-\infty}^x\widehat{N}(\lambda,x,t)q(t)e^{i\lambda\zeta_kt}dt\quad(0\leq k\leq4)\label{eq2.17}
\end{equation}
where
$$\widehat{N}(\lambda,x,t)=\sum\limits_{n=1}^\infty(-i)^n\widehat{K}_n(\lambda,x,t),$$
besides, this series converges in view of estimates \eqref{eq2.16}.

\begin{theorem}\label{t2.2}
When condition \eqref{eq2.2} holds for all $\lambda\in\mathbb{D}_{a/2}$ \eqref{eq2.13}, Jost solutions $\{u_k(\lambda,x)\}_0^2$ of the boundary value problem \eqref{eq2.3}, (b) \eqref{eq2.4} exist and are given by \eqref{eq2.17}. Functions $\{u_k(\lambda,x)\}_0^4$ are holomorphic (in $\lambda$) in the disk $\mathbb{D}_{a/2}$.
\end{theorem}

\begin{remark}\label{r2.2}
Equations \eqref{eq2.5} and \eqref{eq2.14} imply relations
\begin{equation}
\begin{array}{ccc}
{\rm(a)}&v_k(\lambda\zeta_1,x)=v_{k'}(\lambda,x)&(k'=(k+1)(\mod 5));\\
{\rm(b)}&u_k(\lambda\zeta_1,x)=u_{k'}(\lambda,x)&(k'=(k+1)(\mod 5)).
\end{array}\label{eq2.18}
\end{equation}
\end{remark}

Thus, it is sufficient to study properties of the function $v_0(\lambda,x)$ (correspondingly, of $u_0(\lambda,x)$), and properties of other functions $v_p(\lambda,x)$ ($u_p(\lambda,x)$) ($1\leq p\leq4$) are obtained from the properties of $v_0(\lambda x)$ ($u_0(\lambda,x)$) via the substitution $\lambda\rightarrow\lambda\zeta_1^p$.
\vspace{5mm}

{\bf 2.3} Proceed to the analytical (in $\lambda$) properties of the functions $\{v_k(\lambda,x)\}_0^4$. Define the functions
\begin{equation}
\psi_k(\lambda,x)\stackrel{\rm def}{=}v_k(\lambda,x)e^{-i\lambda\zeta_kx}\quad(0\leq k\leq4),\label{eq2.19}
\end{equation}
which, due to \eqref{eq2.5}, are solutions to integral equations
\begin{equation}
\psi_k(\lambda,x)=1+i\int\limits_x^\infty e^{i\lambda\zeta_k(t-x)}\frac{s_4(i\lambda(x-t))}{(i\lambda)^4}q(t)\psi_k(\lambda,t)dt\quad(0\leq k\leq4).\label{eq2.20}
\end{equation}
Hence it follows that
\begin{equation}
\psi_k(\lambda\zeta_1,x)=\psi_{k'}(\lambda,x)\quad(k'=(k+1)(\mod 5),\,0\leq k\leq4)\label{eq2.21}
\end{equation}
due to \eqref{eq2.18}.

Study properties of $\psi_0(\lambda,x)$ which is the solution to equation
\begin{equation}
\psi_0(\lambda,x)=1+i\int\limits_x^\infty e^{i\lambda(t-x)}\frac{s_4(i\lambda(x-t))}{(i\lambda)^4}q(t)\psi_0(\lambda,t)dt.\label{eq2.22}
\end{equation}
Definition of $s_4(z)$ \eqref{eq1.4} implies that
\begin{equation}
\begin{array}{ccc}
{\displaystyle e^{i\lambda(t-x)}\frac{s_4(i\lambda(x-t))}{(i\lambda)^4}=\frac1{5\lambda^4}\left\{1+\frac1{\zeta_1^4}e^{i\lambda(\zeta_0-\zeta_1)(t-x)}+\frac1{\zeta_2^4}e^{i\lambda(\zeta_0-\zeta_2)(t-x)}\right.}\\
{\displaystyle\left.+\frac1{\zeta_3^4}e^{i\lambda(\zeta_0-\zeta_3)(t-x)}+\frac1{\zeta_4^4}e^{\i\lambda(\zeta_0-\zeta_4)(t-x)}\right\}.}
\end{array}\label{eq2.23}
\end{equation}
Consider the first and the last exponents of this sum,
\begin{equation}
e^{i\lambda(\zeta_0-\zeta_1)(t-x)};\quad e^{i\lambda(\zeta_0-\zeta_4)(t-x)}.\label{eq2.24}
\end{equation}
For $\lambda=\alpha+i\beta$ ($\alpha$, $\beta\in\mathbb{R}$), taking into account $\zeta_4=\overline{\zeta_1}$, we obtain
$$i\lambda(\zeta_0-\zeta_1)=(i\alpha-\beta)\left[\left(1-\cos\frac{2\pi}5\right)-i\sin\frac{2\pi}5\right]=-\beta\left(1-\cos\frac{2\pi}5\right)+\alpha\sin\frac{2\pi}5$$
$$+i\left[\alpha\left(1-\cos\frac{2\pi}5\right)+\beta\sin
\frac{2\pi}5\right];$$
$$i\lambda(\zeta_0-\zeta_4)=(i\alpha-\beta)\left[1-\cos\frac{2\pi}5+i\sin\frac{2\pi}5\right]=-\beta\left(1-\cos\frac{2\pi}5\right)-\alpha\sin\frac{2\pi}5$$
$$+i\left[\alpha\left(1-\cos\frac{2\pi}5\right)-\beta\sin\frac{2\pi}5\right],$$
and thus
$$\Re i\lambda(\zeta_0-\zeta_1)=2\sin^2\frac{2\pi}{10}[\alpha K_1-\beta];\quad\Re i\lambda(\zeta_0-\zeta_4)=2\sin^2\frac{2\pi}{10}(-\alpha K_1-\beta)$$
where ${\displaystyle K_1\stackrel{\rm def}{=}\cot\frac{2\pi}{10}}$, therefore, for the exponents \eqref{eq2.24},
$$\left|e^{i\lambda(\zeta_0-\zeta_1)(t-x)}\right|<e^{2\sin^2\frac{2\pi}{10}(\alpha K_1-\beta)(t-x)};\quad \left|e^{i\lambda(\zeta_0-\zeta_4)(t-x)}\right|<e^{2\sin^2\frac{2\pi}{10}(-\alpha K_1-\beta_1)(t-x)}.$$
Thus, for all $\lambda\in K_1$, where

\begin{picture}(200,200)
\put(0,100){\vector(1,0){200}}
\put(100,0){\vector(0,1){200}}
\put(100,100){\vector(1,2){50}}
\put(100,100){\vector(-1,2){50}}
\put(100,100){\vector(-3,1){100}}
\put(100,100){\vector(3,1){100}}
\put(105,170){$K_1$}
\qbezier(116,130)(100,140)(84,130)
\end{picture}

\hspace{20mm} Fig. 2

\begin{equation}
K_1\stackrel{\rm def}{=}\{\lambda\in\mathbb{C}:\alpha k_1<\beta,-\alpha k_1<\beta\}\quad\left(k_1\stackrel{\rm def}{=}\cot\frac{2\pi}{10}\right),\label{eq2.25}
\end{equation}
every exponent \eqref{eq2.24} is less than one in modulo ($t>x>0$) and vanishes when $\lambda\rightarrow\infty$ ($\lambda\in K_1$). Region $K_1$ is a cone symmetric with respect to the ray $il_{\zeta_0}$.

Consider another pair of exponents
\begin{equation}
e^{i\lambda(\zeta_0-\zeta_2)(t-x)};\quad e^{i\lambda(\zeta_0-\zeta_3)(t-x)}\label{eq2.26}
\end{equation}
besides, $\zeta_3=\overline{\zeta}_2$. Analogously to the case of the first two exponents, it is easy to prove that for all $\lambda\in K_2$, where
\begin{equation}
K_2\stackrel{\rm def}{=}\{\lambda\in\mathbb{C}:\alpha k_2<\beta,-\alpha k_2<\beta\}\quad(k_2\stackrel{\rm def}{=}\cot\frac{4\pi}{10});\label{eq2.27}
\end{equation}
every exponent \eqref{eq2.26} is less than one in modulo and vanishes when $\lambda\rightarrow\infty$ ($\lambda\in K_2$).

Monotonous decrease of the function $y=\cot x$ on the interval $(0,\pi)$ implies that $k_2<k_1$, and thus $K_2\supset K_1$ (see Fig. 2).

Since ${\displaystyle k_1=\cot\frac{2\pi}{10}=\tan\left(\frac\pi2-\frac{2\pi}{10}\right)=\tan\frac{3\pi}{10}}$. Define $\theta\in\mathbb{T}$, so that ${\displaystyle\arg\theta=\frac{3\pi}{10}}$. Evidently,
\begin{equation}
K_1\stackrel{\rm def}{=}\Omega_0;\quad\Omega_0\stackrel{\rm def}{=}\theta s_0\label{eq2.28}
\end{equation}
where $s_0$ is given by \eqref{eq1.20}.  And let
\begin{equation}
\Omega_k=\zeta_k\Omega_0\quad(0\leq k\leq4).\label{eq2.29}
\end{equation}

\begin{picture}(200,200)
\put(0,100){\vector(1,0){200}}
\put(100,0){\vector(0,1){200}}
\put(200,67){\line(-3,1){200}}
\put(0,67){\line(3,1){200}}
\put(120,150){$\Omega_0$}
\qbezier(153,120)(100,140)(49,120)
\qbezier(100,23)(155,35)(168,80)
\qbezier(45,120)(30,90)(45,85)
\qbezier(150,120)(168,110)(168,80)
\qbezier(44,80)(40,60)(100,25)
\put(0,105){$\Omega_1$}
\put(140,30){$\Omega_3$}
\put(175,102){$\Omega_4$}
\put(30,27){$\Omega_2$}
\end{picture}

\hspace{20mm} Fig. 3

Sector $\Omega_0$ \eqref{eq2.28} contains 2 rays, $l_{\zeta_1}$ and $\widehat{l}_{\zeta_4}$ \eqref{eq1.18}, besides, $l_{\zeta_1}$ is the bisectrix of the angle between vectors $\theta$ and $i\zeta_0$, and $\widehat{l}_{\zeta_4}$ is also the bisectrix of the angle symmetric to this with regard to the straight line $iL_{\zeta_0}$. Equation \eqref{eq2.9} implies
$$\Omega_1\supset l_{\zeta_2},\widehat{l}_{\zeta_0};\quad\Omega_2\supset l_{\zeta_3},\widehat{l}_{\zeta_1};\quad\Omega_3\supset l_{\zeta_4},\widehat{l}_{\zeta_2};\quad\Omega_4\supset l_{\zeta_0},\widehat{l}_{\zeta_3}.$$
To every function $\{\psi_k(\lambda,x)\}_0^4$ \eqref{eq2.19} there corresponds its holomorphicity sector $\{\Omega_k\}_0^4$ \eqref{eq2.29}. Integral equations \eqref{eq2.20} and analiticity of kernels ${\displaystyle e^{i\lambda(\zeta_k(t-x))}\frac{s_4(i\lambda(x-t))}{(i\lambda)^4}}$ in the corresponding sectors $\{\Omega_k\}_0^4$ \eqref{eq2.29} imply the statement.

\begin{theorem}\label{t2.3}
Functions $\{\psi_k(\lambda,x)\}_0^4$ \eqref{eq2.19} are holomorphic in the corresponding sectors $\{\Omega_k\}_0^4$ \eqref{eq2.29}, besides,
\begin{equation}
\begin{array}{lll}
{\rm(i)}&\lim\limits_{{\lambda\rightarrow\infty}\atop{\lambda\in\Omega_k}}\{\psi_k(\lambda,x)-1\}\cdot\lambda^4=M_5(x);&{\displaystyle M_5(x)\stackrel{\rm def}{=}\frac i{4!}\int\limits_x^\infty(x-t)^4q(t)dt;}\\
{\rm(ii)}&D^5M_5(x)=iq(x)&(\forall x\in\mathbb{R}_+);
\end{array}\label{eq2.30}
\end{equation}
\end{theorem}

\begin{remark}\label{r2.3}
All functions $\{\psi_k(\lambda,x)\}_0^4$ \eqref{eq2.19} are holomorphic in the disk $\mathbb{D}_{a/2}$ \eqref{eq2.13}, but each one has its analytic extension into its (corresponding) sector $\{\Omega_k\}_0^4$. Also note that function $\psi_k(\lambda,x)$ in other sectors $\Omega_s$ ($s\not=k$) does not make sense because of the tendency (growth) to $\infty$ of the kernel of the corresponding integral equation \eqref{eq2.20}.
\end{remark}
\vspace{5mm}

{\bf 2.4} Analogously to \eqref{eq2.19}, define the functions
\begin{equation}
\varphi_k(\lambda,x)=u_k(\lambda,x)e^{-i\lambda\zeta_kx}\quad(0\leq k\leq4).\label{eq2.31}
\end{equation}
By $\Omega_k^-$, we denote the sectors obtained from $\Omega_k$ \eqref{eq2.29} as a result of central symmetry,
\begin{equation}
\Omega_k^-=\{\lambda\in\mathbb{C}:-\lambda\in\Omega_k\}\quad(0\leq k\leq4).\label{eq2.32}
\end{equation}
Analogue of Theorem \ref{t2.3} for the functions $\{\varphi_k(\lambda,x)\}_0^4$ \eqref{eq2.31} looks as follows.

\begin{theorem}\label{t2.4}
Functions $\{\varphi_k(\lambda,x)\}_0^4$ \eqref{eq2.31} are holomorphic (with respect to $\lambda$) at the corresponding sectors $\{\Omega_k^-\}_0^4$ and
\begin{equation}
\begin{array}{lll}
{\rm(i)}&\lim\limits_{{\lambda\rightarrow\infty}\atop{\lambda\in\Omega_k^-}}\{\varphi_k(\lambda,x)-1\}=\widehat{M}_5(x);&{\displaystyle\widehat{M}_5(x)\stackrel{\rm def}{=}-\frac i{4!}\int\limits_{-\infty}^x(x-t)^4q(t)dt;}\\
{\rm(ii)}&D^5\widehat{M}_5(x)=-iq(x)&(\forall x\in\mathbb{R}_-).
\end{array}\label{eq2.33}
\end{equation}
\end{theorem}

For the functions $\{\varphi_k(\lambda,x)\}_0^4$, the following properties similar to the properties of $\{\psi_k(\lambda,x)\}_0^4$, which are given in Remark \ref{r2.3}, are true.

\begin{conclusion}\label{c1}
Knowing any pair of functions $\{\varphi_k(\lambda,x),\varphi_s(\lambda,x)\}$ ($0\leq k$, $s\leq4$) where $\psi_k(\lambda,x)$, as $\varphi_s(\lambda,x)$ also, is one of the elements of the set $\{\psi_k(\lambda,x)\}_0^4$ and $\{\varphi_k(\lambda,x)\}_0^4$ correspondingly, using formulas \eqref{eq2.30} and \eqref{eq2.33}, we can unambiguously restore potential $q(x)$ on the right $x\in\mathbb{R}_+$ and left $x\in\mathbb{R}_-$ half-axes.
\end{conclusion}
\vspace{5mm}

{\bf 2.5} Functions $\{v_k(\lambda,x)\}_0^4$ (as $\{u_k(\lambda,x)\}_0^4$ also) are linearly independent. Really, if for some $c_k\in\mathbb{C}$ ($0\leq k\leq4$) ${\displaystyle\sum\limits_{k=0}^4c_kv_k(\lambda,x)=0}$, then equations \eqref{eq2.6} imply that ${\displaystyle\sum\limits_{k=0}^4c_ke^{i\lambda\zeta_kx}=0}$ ($\forall x\in\mathbb{R}$, $\lambda\not=0$), which is possible only in the case when $c_k=0$ ($\forall k$).

Another proof of linear independence of $\{v_k(\lambda,x)\}_0^4$, elements of which we will need further, is as follows. Define the {\bf matrix Wronskian}
\begin{equation}
W(v)=\left[
\begin{array}{cccc}
v_0(\lambda,x)&v_1(\lambda,x)&...&v_4(\lambda,x)\\
v'_0(\lambda,x)&v'_1(\lambda,x)&...&v'_4(\lambda,x)\\
...&...&...&...\\
v_0^{(4)}(\lambda,x)&v_1^{(4)}(\lambda,x)&...&v_4^{(4)}(\lambda,x)
\end{array}\right],\label{eq2.34}
\end{equation}
and let $w(v)=\det W(v)$ is the standard Wronskian of functions $\col[v_k(\lambda,x),v'_k(\lambda,x),$ $...,v_k^{(4)}(\lambda,x)]$ ($0\leq k\leq4$). Show that $w(v)$ does not depend on $x$; to do this, it is sufficient to establish that ${\displaystyle\frac d{dx}W(v)=0}$. Rewrite $w(v)$ as
$$W(v)=\sum v_k(\lambda,x)w_k(v)$$
where $w_k(v)$ is algebraic compliments of the elements $v_k(\lambda,x)$. Then it is obvious that
$$W'(v)=\sum\limits_{k=0}^4v_k(\lambda,x)w'_k(v).$$
Repeating these considerations for every minor $w_k(v)$, we arrive at derivatives of the elements of the last row of $W(v)$ \eqref{eq2.34}. And since $v_k^{(5)}(\lambda,x)=i(\lambda^5-q(x)v_k(\lambda,x))$ ($0\leq k\leq4$), due to \eqref{eq2.3}, then we obtain that $W'(v)$ has proportionate rows, and thus $W'(v)=0$. So, $w(v)$ does not depend on $x$, and taking into account asymptotics (a) \eqref{eq2.4}, we have
\begin{equation}
W(v)=(i\lambda)^{10}W(\zeta),\label{eq2.35}
\end{equation}
where $W(\zeta)$ is the Vandermonde determinant,
$$W(\zeta)=\det\left|
\begin{array}{cccc}
1&1&...&1\\
\zeta_0&\zeta_1&...&\zeta_4\\
...&...&...&...\\
\zeta_0^4&\zeta_1^4&...&\zeta_4^4
\end{array}\right|$$
and it equals
$$W(\zeta)=\prod\limits_{1\leq z_+<\zeta\leq4}(\zeta_i-\zeta_j).$$
Using $\zeta_k=\zeta_1^k$, we obtain that
\begin{equation}
W(\zeta)=(\zeta_1-1)^4(\zeta_2-1)^3(\zeta_3-1)^2(\zeta_4-1).\label{eq2.36}
\end{equation}
Hence, in view of \eqref{eq2.35}, we find
\begin{equation}
W(v)=(i\lambda)^{10}(\zeta_1-1)^4(\zeta_2-1)^3(\zeta_3-1)^2(\zeta_4-1).\label{eq2.37}
\end{equation}
So, for all $\lambda\not=0$, $W(v)\not=0$, which gives linear independence of $\{v_k(\lambda,x)\}_0^4$.

\section{Scattering problem for the waves incident from ``$+\infty$''}\label{s3}

{\bf 3.1} Every Jost solution $\{u_k(\lambda,x)\}_0^4$ we expand by the other Jost solutions $\{v_k(\lambda,x)\}_0^4$,
\begin{equation}
u_k(\lambda,x)=\sum\limits_{s=0}^4t_{k,s}(\lambda)v_s(\lambda,x)\quad(0\leq k\leq4),\label{eq3.1}
\end{equation}
or, in the matrix form,
\begin{equation}
u(\lambda,x)=T(\lambda)v(\lambda,x)\label{eq3.2}
\end{equation}
where $u(\lambda,x)=\col[u_0(\lambda,x),...,u_4(\lambda,x)]$; $v(\lambda,x)=\col[v_0(\lambda,x),...,v_4(\lambda,x)]$, and $T(\lambda)$ is the {\bf transition matrix},
\begin{equation}
T(\lambda)\stackrel{\rm def}{=}\left[
\begin{array}{cccc}
t_{0,0}(\lambda)&...&t_{0,4}(\lambda)\\
t_{1,0}(\lambda)&...&t_{1,4}(\lambda)\\
...&...&...\\
t_{4,0}(\lambda)&...&t_{4,4}(\lambda)
\end{array}\right].\label{eq3.3}
\end{equation}

\begin{remark}\label{r3.1}
From all the equalities in \eqref{eq3.1}, in view of \eqref{eq2.8}, we should consider one of them (e.g., corresponding to $k=0$),
\begin{equation}
u_0(\lambda,x)=\sum\limits_{s=0}^4t_{0,s}v_s(\lambda,x),\label{eq3.4}
\end{equation}
and all other equalities in \eqref{eq3.1} follow from \eqref{eq3.4} after the substitutions $\lambda\rightarrow\lambda\zeta_k$, besides,
\begin{equation}
\begin{array}{cccc}
t_{1,0}(\lambda)=t_{0,4}(\lambda\zeta_1);&t_{1,1}(\lambda)=t_{0,1}(\lambda\zeta_1);&...&t_{1,4}(\lambda)=t_{0,3}(\lambda\zeta_1);\\
t_{2,0}(\lambda)=t_{0,3}(\lambda\zeta_2);&t_{2,1}(\lambda)=t_{0,4}(\lambda\zeta_2);&...&t_{2,4}(\lambda)=t_{0,2}(\lambda\zeta_2);\\
...&...&...&...\\
t_{4,0}(\lambda)=t_{0,1}(\lambda\zeta_4);&t_{4,1}(\lambda)=t_{0,2}(\lambda\zeta_4)&...&t_{4,4}(\lambda)=t_{0,0}(\lambda\zeta_4).
\end{array}\label{eq3.5}
\end{equation}
\end{remark}

\begin{definition}\label{d3.1}
The four functions $f(x)$, $g(x)$, $\varphi(x)$, $\psi(x)$ from the class $C^4(\mathbb{R}$, we correspond with the Wronskian
\begin{equation}
W(f,g,\varphi,\psi)\stackrel{\rm def}{=}\det\left[
\begin{array}{cccc}
f(x)&g(x)&\varphi(x)&\psi(x)\\
f'(x)&g'(x)&\varphi'(x)&\psi(x)\\
f''(x)&g''(x)&\varphi''(x)&\psi''(x)\\
f'''(x)&g'''(x)&\varphi'''(x)&\psi'''(x)
\end{array}\right].\label{eq3.6}
\end{equation}
\end{definition}

Wronskian $W(f,g,\varphi,\psi)$ is a measure of linear independence of the functions $f(x)$, $g(x)$, $\varphi(x)$, $\psi(x)$ over the field $\mathbb{C}$ and is a natural generalization of the Wronskian  $W(f,g)=f'(x)g(x)-g'(x)f(x)$ for the case of the $n$th order differential equations ($n=4>2$).

\begin{lemma}\label{l3.1}
Wronskian $W(f,g,\varphi,\psi)$ has the following properties:
\begin{equation}
{\rm(i)}\,W(f,g,\varphi,\psi)=-W(g,f,\varphi,\psi)\Rightarrow W(f,f,\varphi,\psi)=0;\label{eq3.7}
\end{equation}

{\rm(ii)} {\bf homogeneity} $w(\alpha f,g,\varphi,\psi)=\alpha W(f,g,\varphi,\psi)$ ($\forall\alpha\in\mathbb{C}$);;

{\rm(iii)} linearity $W(f+\mu,g,\varphi,\psi)=W(f,g,\varphi,\psi)+W(\mu,g,\varphi,\psi)$ ($\forall f$, $g$, $\mu$, $\varphi$, $\psi\in C^3(\mathbb{R})$.
\end{lemma}

P r o o f. We only need to comment on (iii). Its proof follows from the expansion of the determinant $W(f+\mu,g,\varphi,\psi)$ along the first column using analytic complements. $\blacksquare$

\begin{remark}\label{r3.2}
Equation (i) \eqref{eq3.7} implies that properties of homogeneity (ii) and linearity (iii) hold relative to every argument of the Wronskian $W(f,g,\varphi,\psi)$ \eqref{eq3.6}.
\end{remark}

\begin{definition}\label{d3.2}
By $W_{k,s}(v,\lambda,x)$, we denote the Wronskian equal to the minor obtained from the matrix Wronskian $W(v)$ \eqref{eq2.34} via elimination of $k$th column and $s$th row ($0\leq k$, $s\leq4$).
\end{definition}

Evidently, $W_{k,s}(v,\lambda,x)$ up to $(-1)^{k+s}$ coincides with algebraic complement to the element $v_k^{(s)}(\lambda,x)$ of the matrix $w(v)$ \eqref{eq2.34}.

\begin{lemma}\label{l3.2}
For the elements of $W_{4,s}(v,\lambda,x)$ ($0\leq s\leq4$) corresponding to the last ($4$th) row of the matrix Wronskian $W(v)$ \eqref{eq2.34}, the following equalities hold:
\begin{equation}
\begin{array}{ccc}
{\rm(i)}&W'_{4,s}(v,\lambda,x)=W_{4,s-1}(v,\lambda,x)\quad(1\leq s\leq4);\\
{\rm(ii)}&W'_{4,0}(v,\lambda,x)=-i(\lambda^5-q(\lambda))W_{4,4}(v,\lambda,x).
\end{array}\label{eq3.8}
\end{equation}
\end{lemma}

P r o o f. Rewrite $W_{4,4}(v,\lambda,x)$ as
$$W_{4,4}(v,\lambda,x)=v_0V_0+v_1V_1+v_2V_2+v_3V_3,$$
where $\{V_k\}_0^3$ are the algebraic complements of the elements $\{v_k(\lambda,x)\}_0^3$ of the matrix
$$\left[
\begin{array}{cccc}
v_0&v_1&v_2&v_3\\
v'_0&v'_1&v'_2&v'_3\\
v''_0&v''_1&v''_2&v''_3\\
v'''_0&v'''_1&v'''_2&v'''_3
\end{array}\right].$$
Calculate $V'_0$ (for $V'_k$, $k=1$, 2, 3 calculations are similar),
$$\frac d{dx}V_0=\frac d{dx}\left|
\begin{array}{ccc}
v'_1&v'_2&v'_3\\
v''_1&v''_2&v''_3\\
v'''_1&v'''_2&v'''_3
\end{array}\right|=v'_1\left|
\begin{array}{ccc}
v''_2&v''_3\\
v'''_2&v'''_3
\end{array}\right|'-v'_2\left|
\begin{array}{ccc}
v''_x&v''_3\\
v'''_1&v'''_3
\end{array}\right|'+v'_3\left|
\begin{array}{ccc}
v''_1&v''_2\\
v'''_1&v'''_2
\end{array}\right|'$$
$$=\left|
\begin{array}{ccc}
v'_1&v'_2&v'_3\\
v''_1&v''_2&v''_3\\
v_1^{(4)}&v_2^{(4)}&v_3^{(4)}
\end{array}\right|.$$
Hence it follows that $W'_{4,4}(v,\lambda,x)=W_{4,3}(v,\lambda,x)$. It is proved similarly that $W'_{4,k}(v,\lambda,x)=W_{4,k-1}(v,\lambda,x)$ ($k=1$, 2, 3).

Proceed to the proof of (ii) \eqref{eq3.8}. Use the previous technique that we used when proving that $W'_{4,4}(v,\lambda,x)=W_{4,3}(v,\lambda,x)$ and take into account that ${\displaystyle\frac d{dx}v_k^{(4)}(\lambda,x)=}$ $i(\lambda^5-q(x))v_k(\lambda,x)$, and then move the last row to the place of the first one, which gives us a minus sign at the right side of equality (ii) \eqref{eq3.8}. $\blacksquare$

Equation \eqref{eq3.8} implies that $W_{4,4}(v,\lambda,x)$ is the solution to equation
\begin{equation}
\left(-\frac1iD\right)^5y(x)+q(x)y(x)=\lambda^5y(x)\label{eq3.9}
\end{equation}
which is obtained from \eqref{eq2.3} via the substitution $\lambda\rightarrow\overline{\lambda}$ with further complex conjugation of coefficients.

In order to find the explicit form of function $W_{4,4}(v,\lambda,x)$, use the asymptotics (a) \eqref{eq2.4} of the functions $\{v_k(\lambda,x)\}_0^4$, then we obtain that ${\displaystyle W_{4,4}(v,\lambda,x)\rightarrow}$ ${\displaystyle(i\lambda)^6}$ ${\displaystyle\times\prod\limits_0^3e^{i\lambda\zeta_kx}\omega_1}$ where $\omega_1$ is the Vandermonde determinant,
$$\omega_1\stackrel{\rm def}{=}\det\left[
\begin{array}{cccc}
1&1&1&1\\
\zeta_0&\zeta_1&\zeta_2&\zeta_3\\
\zeta_0^2&\zeta_1^2&\zeta_2^2&\zeta_3^2\\
\zeta_0^3&\zeta_1^3&\zeta_2^3&\zeta_3^3
\end{array}\right]=\zeta_4(\zeta_1-1)^3(\zeta_2-1)^2(\zeta_3-1)$$
and using that ${\displaystyle\sum\limits_0^4\zeta_k=0}$, we finally obtain
$$W_{4,4}(v,\lambda,x)\rightarrow(i\lambda)^6\omega_1e^{-i\lambda\zeta_4x}.$$
Define the ``+'' involution,
\begin{equation}
f^+(\lambda)\stackrel{\rm def}{=}\overline{f(\overline{\lambda})}.\label{eq3.10}
\end{equation}
And since $\zeta_4=\overline{\zeta}_1$, then
$$W_{4,4}(\lambda,x)=(i\lambda)^6\omega_1v_1^+(\lambda,x).$$
Properties of the $k$th column ($k=0$, 1, 2, 3) $\{W_{k,s}(v,\lambda,x)\}_{s=0}^4$ of the matrix Wronskian $W(v)$ \eqref{eq2.34} are proved analogously to Lemma \ref{l3.2}.

\begin{theorem}\label{t3.1}
Wronskians $\{W_{k,s}(v,\lambda,x)\}_{s=0}^4$ of the $k$th column ($0\leq k\leq 4$) of the matrix Wronskian $W(v)$ \eqref{eq2.34} have the following properties:
\begin{equation}
\begin{array}{lll}
{\rm(i)}\,W'_{k,s}(v,\lambda,x)=W_{k,s-1}(v,\lambda,x)\quad(1\leq s\leq4,0\leq k\leq4);\\
{\rm(ii)}\,W'_{k,0}(v,\lambda,x)=-i(\lambda^5-q(x))W_{k,4}(v,\lambda,x)\quad(0\leq k\leq4).
\end{array}\label{eq3.11}
\end{equation}
Functions $W_{k,4}(v,\lambda,x)$ ($0\leq k\leq4$) are solutions to equation \eqref{eq3.9} and are given by
$${\rm(a)}\,W_{4,4}(v,\lambda,x)=(i\lambda)^6\omega_1v_1^+(\lambda,x)\quad(\omega_1\stackrel{\rm def}{=}\zeta_4(\zeta_1-1)^2(\zeta_2-1)^2(\zeta_3-1));$$
$${\rm(b)}\,W_{3,4}(v,\lambda,x)=(i\lambda)^6\omega_2v_2^+(\lambda,x)\quad(\omega_2\stackrel{\rm def}{=}\zeta_4(\zeta_1-1)^2(\zeta_2-1)^2(\zeta_3-1)(\zeta_4-1));$$
\begin{equation}
{\rm(c)}\,W_{2,4}(v,\lambda,x)=(i\lambda)^6\omega_3v_3^+(\lambda,x)\quad(\omega_3\stackrel{\rm def}{=}(\zeta_1-1)^2(\zeta_2-1)(\zeta_3-1)^2(\zeta_4-1));\label{eq3.12}
\end{equation}
$${\rm(d)}\,W_{1,4}(v,\lambda,x)=(i\lambda)^6\omega_4v_4^+(\lambda,x)\quad(\omega_4\stackrel{\rm def}{=}\zeta_2(\zeta_1-1)^2(\zeta_2-1)^2(\zeta_3-1)(\zeta_4-1));$$
$${\rm(e)}\,W_{0,4}(v,\lambda,x)=(i\lambda)^6\omega_0v_0^+(\lambda,x)\quad(\omega_0\stackrel{\rm def}{=}(\zeta_1-1)^3(\zeta_2-1)^2(\zeta_3-1)).$$
\end{theorem}

Analogously to \eqref{eq2.34}, we define the matrix Wronskian corresponding to the Jost solutions $\{u_k(\lambda,x)\}_0^4$,
\begin{equation}
W(u)\stackrel{\rm def}{=}
\left[
\begin{array}{cccc}
u_0&u_1&...&u_4\\
u'_0&u'_1&...&u'_4\\
u''_0&u''_1&...&u''_4\\
u_0^{(4)}&u_1^{(4)}&...&u_4^{(4)}
\end{array}\right]\label{eq3.13}
\end{equation}
By $W_{k,s}(u,\lambda,x)$, we define the Wronskians that are equal to the minors obtained from $W(u)$ \eqref{eq3.13} by elimination of the $k$th column and $s$th row. For the Wronskians $W_{k,s}(u,\lambda,x)$, an analogue to Theorem \ref{t3.1} holds.
\vspace{5mm}

{\bf 3.2} Wronskian $W_{0,4}(u,\lambda,x)$ corresponding to $W(u)$ \eqref{eq3.13}, due to (3.12.e), equals
\begin{equation}
W_{0,4}(u,\lambda,x)=\det\left[
\begin{array}{cccc}
u_1(\lambda,x)&...&u_4(\lambda,x)\\
u'_1(\lambda,x)&...&u'_4(\lambda,x)\\
u''_1(\lambda,x)&...&u''_4(\lambda,x)\\
u'''_1(\lambda,x)&...&u'''_4(\lambda,x)
\end{array}\right]=(i\lambda)^6\omega_0u_0^+(\lambda,x)\label{eq3.14}
\end{equation}
Substitute expressions for $\{u_k(\lambda,x)\}_1^4$ into \eqref{eq3.14}, then, taking into account linearity and homogeneity of determinant \eqref{eq3.14} (Lemma \ref{l3.1}), we have
\begin{equation}
(i\lambda)^6\omega_0u_0^+(\lambda,x)=(i\lambda)^6\sum\limits_{p=0}^4\omega_lv_l^+(\lambda,x)\widehat{T}_{0,l}(\lambda)\label{eq3.15}
\end{equation}
where $\widehat{T}_{k,s}(\lambda)$ are algebraic components of the elements $t_{k,s}(\lambda)$ of the matrix $T(\lambda)$ \eqref{eq3.3}.

\begin{remark}\label{r3.2}
The fact that $\widehat{T}_{0,l}(\lambda)$ are algebraic complements of the elements $t_{0,l}(\lambda)$ is proved using induction on the dimension of matrix $T(\lambda)$.
\end{remark}

Rewrite equality \eqref{eq3.15} as
\begin{equation}
u_0^+(\lambda,x)=\sum\limits_{l=0}^4\mu_lv_l^+(\lambda,x)\widehat{T}_{0,l}(\lambda);\label{eq3.16}
\end{equation}
where $\mu_l=\omega_l/\omega_0\in\mathbb{T}$ ($0\leq l\leq4$), besides,
\begin{equation}
\mu_0\stackrel{\rm def}{=}1;\quad\mu_1\stackrel{\rm def}{=}\zeta_4;\quad\mu_2\stackrel{\rm def}{=}\frac{\zeta_4-1}{\zeta_1-1};\quad\mu_3\stackrel{\rm def}{=}\frac{\zeta_3-1}{\zeta_2-1}\cdot\frac{\zeta_4-1}{\zeta_1-1};\quad\mu_4\stackrel{\rm def}{=}\frac{\zeta_4-1}{\zeta_1-1}.\label{eq3.17}
\end{equation}
The following equalities are proved similarly:
\begin{equation}
\mu_ku_k^+(\lambda,x)=\sum\limits_{l=0}^4\mu_lv'_l(\lambda,x)\widehat{T}_{k,l}(\lambda)\quad(1\leq k\leq4).\label{eq3.18}
\end{equation}
Equalities \eqref{eq3.16}, \eqref{eq3.18} in the matrix form become
\begin{equation}
Ju^+(\lambda,x)=\widehat{T}(\lambda)Jv^+(\lambda,x)\label{eq3.19}
\end{equation}
where $u^+(\lambda,x)=\col[u_0^+(\lambda,x),...,u_4^+(\lambda,x)]$; $v^+(\lambda,x)=\col[v_0^+(\lambda,x),...,v_4^+(\lambda,x)]$; $\widehat{T}(\lambda)$ is the matrix formed by algebraic complements $\{\widehat{T}_{k,s}(\lambda)\}$ of the matrix $T(\lambda)$ \eqref{eq3.3}, and $J$ is unitary matrix,
\begin{equation}
J\stackrel{\rm def}{=}\left[
\begin{array}{cccc}
\mu_0&...&0\\
0&\mu_1&...\\
...&...&...&\\
0&...&\mu_4
\end{array}\right]\quad(J^*J=I).\label{eq3.20}
\end{equation}
Since $\widehat{T}(\lambda)=\Delta(\lambda)(T^t(\lambda))^{-1}$ where $\Delta(\lambda)=\det T(\lambda)$, and $T^t(\lambda)$ is the matrix obtained from $T(\lambda)$ via transposition, then equality \eqref{eq3.19} becomes
$$T^t(\lambda)Ju^+(\lambda,x)=\Delta(\lambda)Jv^+(\lambda,x).$$
Applying to both parts of this equality ``+'' involution \eqref{eq3.10}, we obtain
$$T^*(\overline{\lambda})J^*u(\lambda,x)=\Delta^+(\lambda)J^*v(\lambda,x)\Rightarrow\Delta^+(\lambda)v(\lambda,x)=JT^*(\overline{\lambda})J^*u(\lambda,x).$$
Hence, in view of \eqref{eq3.2}, we have
\begin{equation}
\Delta^+(\lambda)J=T(\lambda)JT^*(\overline{\lambda}).\label{eq3.21}
\end{equation}

\begin{theorem}\label{t3.2}
The transition matrix $T(\lambda)$ satisfies relation \eqref{eq3.21} where $J$ is the unitary matrix \eqref{eq3.20} and its determinant $\Delta(\lambda)=\det T(\lambda)$ has the property
\begin{equation}
\Delta(\lambda)=\eta,\quad\eta^3=\zeta_p\label{eq3.22}
\end{equation}
($\zeta_p$ is one of the roots of \eqref{eq1.2}).
\end{theorem}

P r o o f. Considering determinant in equality \eqref{eq3.21}, we obtain
$$(\Delta^+(\lambda))^5\det J=\Delta(\lambda)\Delta^+(\lambda)\det J,$$
which, in view of $\det J\not=0$, gives us
\begin{equation}
(\Delta^+(\lambda))^5=\Delta(\lambda)\Delta^+(\lambda).\label{eq3.23}
\end{equation}
Applying ``+'' to both parts of this equality, we have $(\Delta(\lambda))^5=\Delta^+(\lambda)\Delta(\lambda)$, subtracting this equality from \eqref{eq3.23}, we obtain that $(\Delta(\lambda)/\Delta^+(\lambda))^5=1$ and thus $\Delta(\lambda)\not=\Delta^+(\lambda)=\zeta_p$ ($\zeta_p$ is one of the roots of \eqref{eq1.2}); therefore, \eqref{eq2.23} implies that $(\Delta^+(\lambda))^2\cdot[(\Delta^+(\lambda))^3-\zeta_p]=0$. Equality $\Delta^+(\lambda)=0$ contradicts linear independence of the functions $\{u_k(\lambda,x)\}_0^4$, therefore, $(\Delta^+(\lambda))^3=\zeta_p$. $\blacksquare$
\vspace{5mm}

{\bf 3.3} Equation \eqref{eq3.4} implies
\begin{equation}
r_0u_0(\lambda,x)=v_0(\lambda,x)+\sum\limits_{l=1}^4s_l(\lambda)v_l(\lambda,x)\label{eq3.24}
\end{equation}
where
\begin{equation}
r_0(\lambda)\stackrel{\rm def}{=}1/t_{0,0}(\lambda);\quad s_l(\lambda)\stackrel{\rm def}{=}t_{0,l}(\lambda)/t_{0,0}(\lambda)\quad(1\leq l\leq4).\label{eq3.25}
\end{equation}

\begin{remark}\label{r3.3}
Function $r_0(\lambda)u_0(\lambda,x)$ \label{eq3.24}, in view of (a) \eqref{eq2.4}, behaves in the following way at ``$+\infty$'':
$$r_0(\lambda)u_0(\lambda,x)\rightarrow e^{i\lambda\zeta_0x}+\sum\limits_{l=1}^4s_l(\lambda)e^{i\lambda\zeta_lx}\quad(x\rightarrow\infty).$$
So, {\bf the incident from $+\infty$ wave} $e^{i\lambda\zeta_0x}$ has the {\bf reflected (scattered) wave} ${\displaystyle\sum\limits_1^4s_l(\lambda)e^{i\lambda\zeta_lx}}$, besides, $s_l(\lambda)$ is the {\bf scattering coefficient} of the wave $e^{i\lambda\zeta_0x}$.

At ``$-\infty$'', function $r_0(\lambda)u_0(\lambda,x)$ \eqref{eq3.24}, due to (2.4.b), has the following asymptotic:
$$r_0(\lambda)u_0(\lambda,x)\rightarrow r_0(\lambda)e^{i\lambda\zeta_0x}\quad(x\rightarrow-\infty).$$
Therefore, it is natural to consider $r_0(\lambda)$ as a {\bf coefficient of transmission} of the wave $e^{i\lambda\zeta_0x}$.

Upon equating $(1,1)$ elements in equality \eqref{eq3.21}, we obtain
$$\overline{\eta}=\sum t_{0,l}(\lambda)\mu_lt_{0,l}^+(\lambda)\quad(\eta\in\mathbb{T},\eta^3=\zeta_p)$$
which, in view of \eqref{eq3.25}, gives us
\begin{equation}
\overline{\eta}r_0(\lambda)r_0^+(\lambda)=1+\sum\limits_1^4s_l(\lambda)\mu_ls^+(\lambda).\label{eq3.26}
\end{equation}
Equality \eqref{eq3.26} is the {\bf unitarity property} of the scattering problem for the incident wave $e^{i\lambda\zeta_0x}$.
\end{remark}

Scattering problem for the wave $e^{i\lambda\zeta_kx}$ incident from ``$+\infty$'' ($1\leq k\leq4$) follows from (3.24) after the substitution $\lambda\rightarrow\lambda\zeta_k$, in view of \eqref{eq2.18} and \eqref{eq3.5},
\begin{equation}
r_k(\lambda)u_k(\lambda,x)=v_k(\lambda,x)+\sum\limits_{l\not=k}s_{k,l}(\lambda)v_l(\lambda,x)\quad(1\leq k\leq4)\label{eq3.27}
\end{equation}
where
$$r_k=1/t_{k,k}(\lambda),\quad s_{k,l}(\lambda)=t_{k,l}(\lambda)/t_{k,k}(\lambda)\quad(l\not=k).$$
Define the matrix functions
\begin{equation}
R(\lambda)\stackrel{\rm def}{=}\left[
\begin{array}{cccc}
r_0(\lambda)&0&...&0\\
0&r_1(\lambda)&...&0\\
...&...&...&...\\
0&0&...&r_4(\lambda)
\end{array}\right];\quad S(\lambda)\stackrel{\rm def}{=}\left[
\begin{array}{cccc}
1&s_1(\lambda)&...&s_4(\lambda)\\
s_4(\lambda\zeta_1)&1&...&s_0(\lambda\zeta_1)\\
...&...&...&...\\
s_1(\lambda\zeta_4)&s_2(\lambda\zeta_4)&...&1
\end{array}\right].\label{eq3.28}
\end{equation}
Since
$$T(\lambda)=R^{-1}(\lambda)S(\lambda),$$
then \eqref{eq3.21}, \eqref{eq3.22} imply
\begin{equation}
S(\lambda)JS^*(\overline{\lambda})=\overline{\eta}R(\lambda)JR^*(\overline{\lambda}).\label{eq3.29}
\end{equation}
This equality is the {\bf conservation law} of the problem of scattering waves $\{e^{i\lambda\zeta_k}\}_0^4$, besides, $R(\lambda)$ is the {\bf transition matrix} and $S(\lambda)$ is the {\bf scattering matrix}.
\vspace{5mm}

{\bf 3.4} Study properties of the function $t_{0,0}(\lambda)$. Equation \eqref{eq3.4} implies the equation system
$$\left\{
\begin{array}{llll}
{\displaystyle u_0(\lambda,x)=\sum\limits_0^4t_{0,l}(\lambda)v_l(\lambda,x);}\\
{\displaystyle u'_0(\lambda,x)=\sum\limits_0^4t_{0,1}(\lambda)v'_l(\lambda,x);}\\
............\\
{\displaystyle u_0^{(4)}(\lambda,x)=\sum\limits_0^4t_{0,l}(\lambda)v_l^{(4)}(\lambda,x),}
\end{array}\right.$$
and thus
\begin{equation}
t_{0,0}(\lambda)=W_0(\lambda,x)/\det W(v),\label{eq3.30}
\end{equation}
where $w(v)$ is given by \eqref{eq2.34} and
$$w_0(\lambda,x)=\det\left[
\begin{array}{cccc}
u_0(\lambda,x)&v_1(\lambda,x)&...&v_4(\lambda,x)\\
u'_0(\lambda,x)&v'_1(\lambda,x)&...&v'_4(\lambda,x)\\
...&...&...&...\\
u_0^4(\lambda,x)&v_1^{(4)}(\lambda,x)&...&v_4^{(4)}(\lambda,x)
\end{array}\right].$$
Determinant $\det W(v)$ is found in Section \ref{s2} (see \eqref{eq2.37}), and
$$W_0(\lambda,x)=u_0(\lambda,x)w_{0,0}(v,\lambda,x)-u'_0(\lambda,x)w_{0,1}(\lambda,x)+u''_0(\lambda,x)w_{0,2}(\lambda,x)$$
$$-u'''_0(\lambda,x)w_{0,3}(v,\lambda,x)+u_0^4(\lambda,x)w_{0,4}(v,\lambda,x).$$
Using \eqref{eq3.11} and (3.12.e), we obtain
$$W_0(\lambda,x)=\omega_0(i\lambda)^6\{u_0(\lambda,x)(v_0^+(\lambda,x))^{(4)}-u'_0(\lambda,x)(v_0^+(\lambda,x))'''+u''_0(\lambda,x)(v_0^+(\lambda,x))''$$
$$-u'''_0(\lambda,x)(v_0^+(\lambda,x))'+u_0^{(4)}(\lambda,x)v_0^+(\lambda,x)\}.$$
Assuming that $x=0$ in this equality and taking into account \eqref{eq3.30} and \eqref{eq2.28}, we have
\begin{equation}
t_{0,0}(\lambda)=\frac{F(\lambda)}{c\lambda^4}\quad(c\stackrel{\rm def}{=}\prod\limits_{k=1}^4(\zeta_k-1)>0)\label{eq3.31}
\end{equation}
where
\begin{equation}
\begin{array}{ccc}
F(\lambda)\stackrel{\rm def}{=}u_0(\lambda,0)(v_0^+(\lambda,0))^{(4)}-u'_0(\lambda,0)(v_0^+(\lambda,0))'''+u_0''(\lambda,0)(v_0^+(\lambda,0))''\\
-u'''_0(\lambda,0)(v_0^+(\lambda,0))'+u_0^{(4)}(\lambda,0)v_0^+(\lambda,0).
\end{array}\label{eq3.32}
\end{equation}
Holomorphicity of $u_0(\lambda,0)$ and $v_0^+(\lambda,0)$ in the sector $\Omega_0^-$ implies holomorphicity of the function $F(\lambda)$ \eqref{eq3.32} in $\Omega_0^-$. Moreover, for $\lambda\rightarrow\infty$ ($\lambda\in\Omega_0^-$), due to \eqref{eq2.5} and \eqref{eq2.12},
$$u_0(\lambda,0)(v_0^+(\lambda,0))^{(4)}\rightarrow\lambda^4;\quad u'_0(\lambda,0)(v_0^+(\lambda,0))'''\rightarrow\lambda^4,\quad u''_0(\lambda,0)(v_0^+(\lambda,0))''\rightarrow\lambda^4;$$
$$u'''_0(\lambda,0)(v_0^+(\lambda,0))'\rightarrow\lambda^4;\quad u_0^{(4)}(\lambda,0)v_0^+(\lambda,0)\rightarrow\lambda^4,$$
therefore, for $t_{0,0}(\lambda)$ the following asymptotic is true:
\begin{equation}
t_{0,0}(\lambda)\rightarrow\frac1c\quad(\lambda\rightarrow\infty,\lambda\in\Omega_0^-).\label{eq3.33}
\end{equation}
Therefore, zeros of $t_{0,0}(\lambda)$ \eqref{eq3.31} lie in the bounded region in $\Omega_0^-$ and there is a finite number of such zeros, due to analyticity of $t_{0,0}(\lambda)$ in $\Omega_0^-$.

To describe the multiplicity of zeros of the function $t_{0,0}(\lambda)$, we define function $A(\lambda)$,
\begin{equation}
\begin{array}{cccc}
{\displaystyle A(\lambda)\stackrel{\rm def}{=}u_0(\lambda,0)(v_0^+(\lambda,0))^{(4)}\frac d{d\lambda}\ln\frac{u_0(\lambda,0)(v_0^+(\lambda,0))^{(4)}}{u''_0(\lambda,0)(v_0^+(\lambda,0))''}}\\
{\displaystyle-u'_0(\lambda,0)(v_0^+(\lambda,0))'''\frac d{d\lambda}\ln\frac{u_0^+(\lambda,0)(v_0^+(\lambda,0))'''}{u''_0(\lambda,0)(v_0^+(\lambda,0))''}}\\
{\displaystyle-u'''_0(\lambda,0)(v_0^+(\lambda,0))'\frac d{d\lambda}\ln\frac{u'''_0(\lambda,0)(v_0^+(\lambda,0))'}{u''_0(\lambda,0)(v_0^+(\lambda,0))''}}\\
{\displaystyle+u_0^{(4)}(\lambda,0)v_0^+(\lambda,0)\frac d{d\lambda}\ln\frac{u_0^{(4)}(\lambda,0)v_0^+(\lambda,0)}{u''_0(\lambda,0)(v_0^+(\lambda,0))''}}.
\end{array}\label{eq3.34}
\end{equation}

\begin{lemma}\label{l3.3}
Zero $w$ of the function $F(\lambda)$ \eqref{eq3.32} is simple if $A(w)\not=0$. When $F(w)=0$ and $A(w)=0$, zero $w$ of the function $F(\lambda)$ is of multiplicity not less than two.
\end{lemma}

P r o o f. Statement follows from the equality ${\displaystyle\left.\frac d{d\lambda}F(\lambda)\right|_{\lambda=w}=A(w)}$. Since
$$\frac d{d\lambda}F(\lambda)=u_0(\lambda,0)(v_0^+(\lambda,0))^{(4)}\frac d{d\lambda}\ln u_0(\lambda,0)(v_0^+(\lambda,0))^{(4)}$$
$$-u'_0(\lambda,0)(v_0^+(\lambda,0))'''\cdot\frac d{d\lambda}\ln u'_0(\lambda)(v_0^+(\lambda,0))'''$$
$$+u''_0(\lambda,0)(v_0^+(\lambda,0))''\frac d{d\lambda}\ln u''_0(\lambda,0)(v_0^+(\lambda,0))''$$
$$-u'''_0(\lambda,0)(v_0^+(\lambda,0))'\frac d{d\lambda}u'''_0(\lambda,0)(v_0^+(\lambda,0))'$$
$$+u_0^{(4)}(\lambda,0)v_0^+(\lambda,0)\cdot\frac d{d\lambda}\ln u_0^{(4)}(\lambda,0)v_0^+(\lambda,0),$$
then, assuming in this equality $\lambda=w$ and taking into account relation
$$u''_0(w,0)(v_0^+(w,0))''=-u_0(w,0)(v_0^+(w,0))^{(4)}+u'_0(w,0)(v_0^+(w,0))'''$$
$$+u'''_0(w,0)(v_0^+(w,0))'-u_0^{(4)}(w,0)v_0^+(w,0)$$
that follows from $F(w)=0$, we obtain
$$\left.\frac d{d\lambda}F(\lambda)\right|_{\lambda=w}=A(w).\blacksquare$$

\begin{remark}\label{r3.4}
Higher multiplicity ($n>2$) of zeros of the function $F(\lambda)$ \eqref{eq3.32} is studied analogously. Hereinafter, we assume that all zeros of $F(\lambda)$ are simple.
\end{remark}

\begin{remark}\label{r3.5}
If $w\in\Omega_0^-\cap\Omega_2$ is a general zero of the functions $t_{0,0}(\lambda)$, $t_{0,1}(\lambda)$, $t_{0,3}(\lambda)$, $t_{0,4}(\lambda)$ ($t_{0,0}(w)=t_{0,1}(w)=t_{0,3}(w)=t_{-0,4}(w)=0$), then \eqref{eq3.4} implies $u(w,x)=t_{0,2}(w)v_2(w,x)$. Since $u_0(w,x)\in L^2(\mathbb{R}_-)$ and $v_2(w,x)\in L^2(\mathbb{R}_+)$, then $u_0(w,x)\in L^2(\mathbb{R})$ and it is an eigenfunction of operator $L_q$, besides, $w^5$ is real ($L_q$ is self-adjoint). And thus $w$ belongs to the one of rays, either to $\{l_{\zeta_k}\}_0^4$ or $\{\widehat{l}_{\zeta_k}\}_0^4$ \eqref{eq1.12}. And since $\widehat{l}_{\zeta_1}$ is the only ray lying in $\Omega_0^-\cap\Omega_2$, i.e., $w\in\widehat{l}_{\zeta_1}$, and thus $w$ belongs to the set
\begin{equation}
\Lambda_1^-\stackrel{\rm def}{=}\{\widehat{\varkappa}_m\stackrel{\rm def}{=}\zeta_1b_m:b_m<0,\,1\leq m\leq M<\infty\}.\label{eq3.35}
\end{equation}
Analogously, if $w\in\Omega_0^-\cap\Omega_3$ is a joint zero of the functions $t_{0,0}(\lambda)$, $t_{0,1}(\lambda)$, $t_{0,2}(\lambda)$, $t_{0,4}(\lambda)$ ($t_{0,0}(w)=t_{0,1}(w)=t_{0,2}(w)=t_{0,4}(w)=0$), then $u_0(w,x)=t_{0,3}(w)v_3(w,x)$, and thus $u_0(w,x)\in L^2(\mathbb{R})$, i.e., $u_0(w,x)$ is also an eigenfunction of operator $L_q$. Reality of $w^5$ implies that $w\in l_{\zeta_4}$ ($l_{\zeta_4}\in\Omega_0^-\cap\Omega_3$), i.e., $w$ belongs to the set
\begin{equation}
\Lambda_4^+\stackrel{\rm def}{=}\{\varkappa_n\stackrel{\rm def}{=}\zeta_4a_n:a_n>0;\,1\leq n\leq N<\infty\}.\label{eq3.36}
\end{equation}

\begin{picture}(200,200)
\put(0,100){\vector(1,0){200}}
\put(100,0){\vector(0,1){200}}
\put(100,100){\vector(-1,-2){50}}
\put(100,100){\vector(1,-2){50}}
\put(100,100){\vector(3,-1){100}}
\put(100,100){\vector(-3,-1){100}}
\put(76,55){$\circ$}
\put(50,52){$\zeta_1b_m$}
\put(37,7){$\widehat{l}_{\zeta_1}$}
\put(150,7){$l_{\zeta_4}$}
\put(105,30){$\Omega_0^-$}
\put(119,50){$\times$}
\put(128,55){$\zeta_4a_n$}
\qbezier(153,80)(100,60)(49,80)
\end{picture}

\hspace{20mm} {\rm Fig. 4. Zeros of function $t_{0,0}(\lambda)$}

So, set of points in $\Omega_0^-$ corresponding to singly connected states is given by $\Lambda_1^-\cup\Lambda_4^+$ (see Fig. 4).
\end{remark}

\begin{remark}\label{r3.6}
Equation $t_{k,k}(\lambda)=t_{0,0}(\lambda\zeta_k)$ ($0\leq k\leq4$) implies that set of zeros of $t_{k,k}(\lambda)$ is situated in the sector $\Omega_k^-$ and is given by $\zeta_k(\Lambda_1^-\cup\Lambda_4^+$). Singly connected states of operator $L_q$ also correspond to these zeros.
\end{remark}

Equation \eqref{eq3.24}, upon substitutions $\lambda\rightarrow\lambda\zeta_k$, due to \eqref{eq2.18}, imply the 4 additional relations (see also \eqref{eq3.27}) corresponding to the scattering problems for the waves $\{e^{i\lambda\zeta_kx}\}_0^4$ incident from ``$+\infty$'',
$${\rm(i)}\,r_1(\lambda)u_1(\lambda,x)=s_4(\lambda\zeta_1)v_0(\lambda,x)+v_1(\lambda,x)+s_1(\lambda\zeta_1)v_2(\lambda,x)+s_2(\lambda\zeta_1)v_3(\lambda,x)$$
$$+s_3(\lambda\zeta_1)v_4(\lambda,x);$$
\begin{equation}
\begin{array}{cccc}
{\rm(ii)}\,r_2(\lambda)u_2(\lambda,x)=s_3(\lambda\zeta_2)v_0(\lambda,x)+s_4(\lambda\zeta_1)v_1(\lambda,x)+v_2(\lambda,x)+s_1(\lambda\zeta_2)v_3(\lambda,x)\\
+s_2(\lambda\zeta_2)v_4(\lambda,x);\\
{\rm(iii)}\,r_3(\lambda)u_3(\lambda,x)=s_2(\lambda\zeta_3)v_0(\lambda,x)+s_3(\lambda\zeta_3)v_1(\lambda,x)+s_4(\lambda\zeta_3)v_2(\lambda,x)+v_3(\lambda,x)\\
+s_1(\lambda\zeta_3)v_4(\lambda,x);
\end{array}\label{eq3.37}
\end{equation}
$${\rm(iv)}\,r_4(\lambda)u_4(\lambda,x)=s_1(\lambda\zeta_4)v_0(\lambda,x)+s_2(\lambda\zeta_4)v_1(\lambda,x)+s_3(\lambda\zeta_4)v_2(\lambda,x)+s_4(\lambda\zeta_4)v_3(\lambda,x)$$
$$+v_4(\lambda,x)$$
where $r_k(\lambda)\stackrel{\rm def}{=}r_0(\lambda\zeta_k)$ ($0\leq k\leq4$).
\vspace{5mm}

{\bf 3.5} Into the Wronskians $W(u_0,v_1,v_2,v_3)$, $W(u_0,v_2,v_3,v_4)$, we substitute expression \eqref{eq3.24} for the function $u_0(\lambda,x)$, then, using Lemma \ref{l3.1}, we have
$${\rm(i)}\,r_0(\lambda)W(u_0,v_1,v_2,v_3)=W(v_0,v_1,v_2,v_3)-s_4(\lambda)W(v_1,v_2,v_3,v_4);$$
$${\rm(ii)}\,r_0(\lambda)W(u_0,v_2,v_3,v_4)=W(v_0,v_2,v_3,v_4)+s_1(\lambda)W(v_1,v_2,v_3,v_4).$$
Hence, due to (a) \eqref{eq3.12}, (d) \eqref{eq3.12}, (e) \eqref{eq3.12}, we obtain
$${\rm(i)}\,r_0(\lambda)W(u_0,v_1,v_2,v_3)=(i\lambda)^6\{\omega_1v_1^+(\lambda,x)-s_4(\lambda)\omega_0v_0^+(\lambda,x)\};$$
$${\rm(ii)}\,r_0(\lambda)W(u_0,v_2,v_3,v_4)=(i\lambda)^6\{\omega_4v_4^+(\lambda,x)+s_1(\lambda)\omega_0v_0^+(\lambda,x)\}.$$
Taking into account $v_k(\lambda,x)=\psi_k(\lambda,x)e^{i\lambda\zeta_kx}$, we arrive at the relations
\begin{equation}
\begin{array}{ccc}
{\rm(i)}\,f_1(\lambda,x)\omega_0-\omega_1\psi_1^+(\lambda,x)=-\omega_0s_4(\lambda)e^{i\lambda(\zeta_4-\zeta_0)x}\psi_0^+(\lambda,x)&(\widehat{l}_\theta=l_{\overline{\theta}\zeta_1});\\
{\rm(ii)}\,f_4(\lambda,x)\omega_0-\omega_4\psi_4^+(\lambda,x)=\omega_0s_1(\lambda)e^{i\lambda(\zeta_1-\zeta_0)x}\psi_0^+(\lambda,x)&(\widehat{l}_{\theta\zeta_1}=l_{\overline{\theta}}),
\end{array}\label{eq3.38}
\end{equation}
where
\begin{equation}
\begin{array}{ccc}
f_1(\lambda,x)\stackrel{\rm def}{=}r_0(\lambda)W(u_0,v_1,v_2,v_3)e^{i\lambda\zeta_4x}(i\lambda)^{-6};\\
f_4(\lambda,x)\stackrel{\rm def}{=}r_0(\lambda)W(u_0,v_2,v_3,v_4)e^{i\lambda\zeta_1x}(i\lambda)^{-6}.
\end{array}\label{eq3.39}
\end{equation}

\begin{remark}\label{r3.7}
Lemma \ref{l3.1} implies that
$$W(u_0,v_1,v_2,v_3)=W(u_0,v_1+v_2,v_2,v_3+v_2).$$
Every summand of the sum $v_1+v_2$ exists and is analytical in the disk $\mathbb{D}_{a/2}$ (Remark \ref{r2.3}) and it can be holomorphically extended into the sector $\Omega_1$ and $\Omega_2$ correspondingly. Similarly, summands of the sum $v_3+v_2$ are analytical in the sectors $\Omega_3$ and $\Omega_2$. Hence it follows that $W(u_0,v_1,v_2,v_3)$ is holomorphic in the sector $\Omega_0^-\cap\Omega_2$ excluding poles at the points of $\Lambda_1^-$ \eqref{eq3.35}. Holomorphicity of $W(u_0,v_2,v_3,v_4)$ in the sector $\Omega_0^-\cap\Omega_3$ is proved analogously, excluding poles at the points of $\Lambda_4^+$ \eqref{eq3.36}. Equation \eqref{eq3.33} implies asymptotics
\begin{equation}
f_1(\lambda,x)\rightarrow c\omega_0\,(\lambda\rightarrow\infty,\lambda\in\Omega_0^-\cap\Omega_2);\quad f_4(\lambda,x)\rightarrow c\omega_0\,(\lambda\rightarrow\infty,\lambda\in\Omega_0^-\cap\Omega_3).\label{eq3.40}
\end{equation}
\end{remark}

\begin{picture}(200,200)
\put(0,100){\vector(1,0){200}}
\put(100,0){\vector(0,1){200}}
\put(100,100){\vector(-1,-2){50}}
\put(100,100){\vector(1,-2){50}}
\put(100,100){\vector(3,1){100}}
\put(100,100){\vector(-3,1){100}}
\put(50,72){$\psi_1^+$}
\put(85,67){$f_1$}
\put(102,67){$f_4$}
\put(37,7){$\widehat{l}_{\theta}$}
\put(150,7){$\widehat{l}_{\theta\zeta_1}$}
\put(150,50){$\Omega_1^-$}
\put(175,135){$\widehat{l}_{\theta\zeta_2}$}
\put(0,135){$\widehat{l}_{\theta\zeta_4}$}
\put(106,5){$_{\Omega_0^-\cap\Omega_3}$}
\put(60,5){$_{\Omega_0^-\cap\Omega_2}$}
\put(128,75){$\psi_4^+$}
\qbezier(100,50)(80,56)(75,55)
\qbezier(75,55)(20,50)(35,122)
\qbezier(100,50)(110,45)(125,50)
\qbezier(125,50)(160,70)(165,122)
\end{picture}

\hspace{20mm} {\rm Fig. 5. The Riemann problem \eqref{eq3.39}}

Equality (i) \eqref{eq3.38} is a Riemann boundary value problem on the ray $\widehat{l}_{\theta}$ for the functions $\psi_1^+(\lambda,x)$ and $f_1(\lambda,x)$ holomorphic in the sectors $\Omega_4^-$ and $\Omega_0^-\cap\Omega_2$. Similarly, relation (ii) \eqref{eq3.38} is a Riemann boundary value problem on the ray $\widehat{l}_{\theta\zeta_1}$ for the functions $\psi_4^+(\lambda,x)$ and $f_4(\lambda,x)$ holomorphic in the sectors $\Omega_1^-$ and $\Omega_0^-\cap\Omega_3$.
\vspace{5mm}

{\bf 3.6}
\begin{definition}\label{d3.3}
A holomorphic in the sector $\Omega_0$ function $F(\lambda)$ is said to be of class $SYMM(\Omega_0)$ if its values on the boundaries of this sector $\Omega_0$ coincide,
$$F(\theta t)=F(\theta\zeta_1t)\quad(\forall t\in\mathbb{R}_+).$$
\end{definition}

If $v_0(\lambda,x)\in SYMM(\Omega_0)$, then $v_k(\lambda,x)\in SYMM(\Omega_k)$ ($1\leq k\leq4$), in view of \eqref{eq2.18}.

Substitute expression \eqref{eq3.24} for the function $u_0(\lambda,x)$ into the Wronskians $W(u_0,v_1,v_2,v_4)$, $W(u_0,v_1,v_3,v_4)$, then, in accordance with Lemma \ref{l3.1}, we obtain
$${\rm(iii)}\,r_0(\lambda)W(u_0,v_1,v_2,v_4)=W(v_0,v_1,v_2,v_4)+s_3(\lambda)W(v_1,v_2,v_3,v_4);$$
$${\rm(iv)}\,r_0(\lambda)W(u_0,v_1,v_3,v_4)=W(v_0,v_1,v_3,v_4)-s_2(\lambda)W(v_1,v_2,v_3,v_4);$$
and taking into account (b) \eqref{eq3.12}, (c) \eqref{eq3.12}, (e) \eqref{eq3.12}, we have
$${\rm(iii)}\,r_0(\lambda)W(u_0,v_1,v_2,v_4)=(i\lambda)^6\{\omega_2v_2^+(\lambda,x)+s_3(\lambda)\omega_0v_0^+(\lambda,x)\};$$
$${\rm(iv)}\,r_0(\lambda)W(u_0,v_1,v_3,v_4)=(i\lambda)^6\{\omega_3v_3^+(\lambda,x)-s_2(\lambda)\omega_0v_0^+(\lambda,x)\}.$$

Assuming that $v_0(\lambda,x)\in SYMM(\Omega_0)$, i.e., $v_0(\theta t,x)=v_0(\theta\zeta_1t,x)$ ($\forall t\in\mathbb{R}_+$). The last statement means that $\left.v_2(\lambda,x)\right|_{l_{\theta\zeta_3}}=\left.v_3(\lambda,x)\right|_{l_{\theta\zeta_3}}$, therefore, the left parts of equalities (iii) and (iv) coincide for all $\lambda\in l_{\theta\zeta_3}$. Consequently, upon subtraction, we obtain
\begin{equation}
\omega_3v_3^+(\lambda,x)-\omega_2v_2^+(\lambda,x)=[s_2(\lambda)+s_3(\lambda)]\omega_0v_0^+(\lambda,x)\quad(\forall\lambda\in\widehat{l}_{\theta\zeta_3}).\label{eq3.41}
\end{equation}
As a result, we arrive at the Riemann problem on the ray $\widehat{l}_{\theta\zeta_3}$ for the functions $v_3^+(\lambda,x)$ and $v_2^+(\lambda,x)$ holomorphic in the sectors $\Omega_2^-$ and $\Omega_3^-$ (Fig.6).

\begin{picture}(200,200)
\put(0,100){\vector(1,0){200}}
\put(100,0){\vector(0,1){200}}
\put(100,100){\vector(-1,-2){50}}
\put(100,100){\vector(1,-2){50}}
\put(100,100){\vector(3,1){100}}
\put(100,100){\vector(-3,1){100}}
\put(50,72){$v_1^+$}
\put(150,50){$\Omega_1^-$}
\put(30,50){$\Omega_4^-$}
\put(30,150){$\Omega_3^-$}
\put(130,150){$\Omega_2^-$}
\put(105,185){$\widehat{l}_{\theta\zeta_3}$}
\put(175,135){$\widehat{l}_{\theta\zeta_2}$}
\put(0,135){$\widehat{l}_{\theta\zeta_4}$}
\put(128,75){$v_4^+$}
\qbezier(100,150)(80,156)(35,122)
\qbezier(75,55)(20,50)(35,122)
\qbezier(100,150)(110,155)(165,122)
\qbezier(125,50)(160,70)(165,122)
\end{picture}

\hspace{20mm} {\rm Fig. 6. Riemann problems \eqref{eq3.38}}

Note that substitution $\lambda\rightarrow\lambda\zeta_1$ transforms the triplet $\{v_2^+(\lambda,x),\widehat{l}_{\theta\zeta_3},v_3^+(\lambda,x)\}$ from boundary value problem \eqref{eq3.41} into the triplet $\{v_1^+(\lambda,x),\widehat{l}_{\theta\zeta_4},v_2^+(\lambda,x)\}$, which leads us to the following boundary value problem on the ray $\widehat{l}_{\theta\zeta_4}$:
\begin{equation}
\omega_3v_2^+(\lambda,x)-\omega_2v_1^+(\lambda,x)=[s_2(\lambda\zeta_1)+s_3(\lambda\zeta_1)]\omega_0v_4^+(\lambda,x)\quad(\forall\lambda\in\widehat{l}_{\theta\zeta_4}).\label{eq3.42}
\end{equation}
Analogously, substitution $\lambda\rightarrow\lambda/\zeta_1=\lambda\zeta_4$ transforms triplet $\{v_2(\lambda,x),\widehat{l}_{\theta\zeta_3},v_3^+(\lambda,x)\}$ into the triplet $\{v_3^+(\lambda,x),\widehat{l}_{\theta\zeta_2},v_4^+(\lambda,x)\}$, which gives us the boundary value problem on the ray $\widehat{l}_{\theta\zeta_2}$,
\begin{equation}
\omega_3v_4^+(\lambda,x)-\omega_2v_3^+(\lambda,x)=[s_2(\lambda\zeta_4)+s_3(\lambda\zeta_4)]\omega_0v_1^+(\lambda,x)\quad(\lambda\in\widehat{l}_{\theta\zeta_2}).\label{eq3.43}
\end{equation}
Define the piecewise holomorphic function
\begin{equation}
\mathcal{F}(\psi,\lambda,x)\stackrel{\rm def}{=}\left\{
\begin{array}{cccccc}
f_4(\lambda,x)/\omega_0c&(\lambda\in\Omega_0^-\cap\Omega_3);\\
f_1(\lambda,x)/\omega_0c&(\lambda\in\Omega_0^-\cap\Omega_2);\\
\psi_1^+(\lambda,x)&(\lambda\in\Omega_4^-);\\
\psi_2^+(\lambda,x)&(\lambda\in\Omega_3^-);\\
\psi_3^+(\lambda,x)&(\lambda\in\Omega_2^-);\\
\psi_4^+(\lambda,x)&(\lambda\in\Omega_1^-).
\end{array}\right.\label{eq3.44}
\end{equation}
Function $\mathcal{F}(\psi,\lambda,x)$ \eqref{eq3.44} is holomorphic in the sectors $\{\Omega_k^-\}_1^4$ and meromorphic in $\Omega_0^-\cap\Omega_2$, $\Omega_0^-\cap\Omega_3$. Moreover, $\mathcal{F}(\psi,f,\overline{\lambda})\rightarrow1$ as $\lambda\rightarrow\infty$ (inside of every sector).

Write relations (i) \eqref{eq3.38}, (ii) \eqref{eq3.38}, and \eqref{eq3.41} -- \eqref{eq3.43} as jump problems
$${\rm(a)}\,\psi_1^+(\lambda,x)-f_1(\lambda,x)=s_4(\lambda)e^{i\lambda(\zeta_4-\zeta_0)x}\psi_0^+(\lambda,x)+\left(1-\frac{\omega_1}{\omega_0}\right)\psi_1^+(\lambda,x)\quad(\widehat{l}_\theta);$$
$${\rm(b)}\,\psi_4^+(\lambda,x)-f_4(\lambda,x)=-s_1(\lambda)e^{i\lambda(\zeta_1-\zeta_0)x}\psi_0^+(\lambda,x)+\left(1-\frac{\omega_4}{\omega_0}\right)\psi_4^+(\lambda,x)\quad(\widehat{l}_{\theta\zeta_1});$$
\begin{equation}
\begin{array}{ccc}
{\displaystyle{\rm(c)}\,\psi_3^+(\lambda,x)-\psi_2^+(\lambda,x)=\frac{\omega_0}{\omega_3}[s_2(\lambda)+s_3(\lambda)]e^{i\lambda(\zeta_2-\zeta_0)x}\psi_0^+(\lambda,x)}\\
{\displaystyle+\left(\frac{\omega_2}{\omega_3}e^{i\lambda(\zeta_2-\zeta_3)x}-1\right)\psi_2^+(\lambda,
x)\quad(\widehat{l}_{\theta\zeta_3});}
\end{array}\label{eq3.45}
\end{equation}
$${\rm(d)}\,\psi_2^+(\lambda,x)-\psi_1^+(\lambda,x)=\frac{\omega_0}{\omega_3}[s_2(\lambda\zeta_1)+s_3(\lambda\zeta_1)]e^{i\lambda(\zeta_3-\zeta_1)}\psi_4^+(\lambda,x)$$
$$+\left(\frac{\omega_2}{\omega_3}e^{i\lambda(\zeta_3-\zeta_4)x}-1\right)\psi_1^+(
\lambda,x)\quad(\widehat{l}_{\theta\zeta_4});$$
$${\rm(e)}\,\psi_4^+(\lambda,x)-\psi_3^+(\lambda,x)=\frac{\omega_0}{\omega_3}[s_2(\lambda\zeta_4)+s_3(\lambda\zeta_4)]e^{i\lambda(\zeta_1-\zeta_4)x}\psi_1^+(\lambda,x)$$
$$+\left(\frac{\omega_2}{\omega_3}e^{i\lambda(\zeta_1-\zeta_2)x}-1
\right)\psi_3^+(\lambda,x)\quad(\widehat{l}_{\theta\zeta_2}).$$
Explicit expression of the function $\mathcal{F}(\psi,\lambda,x)$ \eqref{eq3.44} follows from the Cauchy formula
$$\mathcal{F}(\psi,\lambda,x)=1+\frac1{2\pi i}\int\limits_\gamma\frac{\mathcal{F}(\psi,w,x)-1}{w-\lambda}dw$$
where $\gamma$ is the contour depicted on Fig. 7.

\begin{picture}(400,380)
\put(0,200){\vector(1,0){400}}
\put(197,205){\vector(0,1){50}}
\put(197,255){\line(0,1){50}}
\put(203,305){\vector(0,-1){50}}
\put(203,255){\line(0,-1){50}}
\put(200,0){\vector(0,1){380}}
\put(188,200){\vector(-3,1){88}}
\put(100,230){\line(-3,1){30}}
\put(66,250){\vector(3,-1){35}}
\put(100,239){\line(3,-1){100}}
\put(188,200){\vector(-3,-1){35}}

\put(200,206){\vector(3,1){100}}
\put(300,240){\line(3,1){44}}
\put(344,242){\vector(-3,-1){42}}
\put(300,226){\line(-3,-1){70}}
\put(200,200){\line(-3,1){200}}
\put(200,200){\line(3,1){200}}
\put(4,280){$\widehat{l}_{\theta\zeta_4}$}
\put(202,310){$\widehat{l}_{\zeta_3}$}
\qbezier(344,256)(200,360)(66,250)
\qbezier(344,242)(200,0)(68,242)
\put(360,240){$\widehat{l}_{\theta\zeta_2}$}
\put(120,120){$\Omega_0^-\cap\Omega_3$}
\put(270,120){$\Omega_0^-\cap\Omega_3$}
\put(179,166){$\times$}
\put(210,166){$\times$}
\put(227,130){$\times$}
\put(161,130){$\times$}
\put(103,185){$\psi_1^+$}
\put(123,185){$\Omega_4^-$}
\put(90,170){$\widehat{l}_\theta$}
\put(150,280){$\psi_2^+\Omega_3^-$}
\put(250,280){$\psi_3^+$}
\put(265,265){$\Omega_2^-$}
\put(330,280){$\gamma$}
\put(80,205){$\gamma$}
\put(230,200){\vector(3,-1){30}}
\put(250,193){\line(3,-1){45}}
\put(200,200){\line(3,-1){90}}
\put(200,200){\line(-3,-1){80}}
\put(154,190){\line(-3,-1){45}}
\put(200,195){\vector(-3,-1){40}}
\put(162,182){\line(-3,-1){40}}
\put(200,195){\vector(3,-1){40}}
\put(238,182){\line(3,-1){40}}
\put(205,194){\vector(1,-2){11}}
\qbezier(216,174)(223,171)(219,167)
\put(219,167){\line(1,-2){5}}
\multiput(222,155)(2,-4){3}{.}
\put(228,148){\vector(1,-2){5}}
\qbezier(233,138)(243,134)(235,130)
\qbezier(235,130)(223,130)(229,137)
\put(229,137){\line(-1,2){5}}
\multiput(222,147)(-2,4){3}{.}
\put(220,155){\vector(-1,2){5}}
\qbezier(214,164)(205,165)(210,174)
\put(210,174){\line(-1,2){10}}
\put(200,194){\line(-1,-2){10}}
\qbezier(190,175)(195,170)(186,165)
\put(186,165){\line(-1,-2){5}}
\multiput(177,153)(-2,-4){3}{.}
\put(175,145){\vector(-1,-2){5}}
\qbezier(171,136)(175,133)(166,130)
\qbezier(166,130)(155,135)(167,138)
\put(167,138){\line(1,2){5}}
\multiput(169,148)(2,4){3}{.}
\put(174,156){\vector(1,2){5}}
\qbezier(180,166)(173,170)(185,175)
\put(185,175){\line(1,2){11}}
\put(283,150){$\widehat{l}_{\theta\zeta_1}$}
\put(273,203){$\psi_4^+$}
\put(295,203){$\Omega_1^-$}
\put(230,140){$f_4/\omega_0c$}
\put(175,130){$f_1/\omega_0c$}
\end{picture}

\hspace{20mm} Fig. 7. Contour of the Cauchy integral

Using \eqref{eq3.45}, we obtain the representation of $\mathcal{F}(\psi,\lambda,x)$ \eqref{eq3.44},
\begin{equation}
\mathcal{F}(\psi,\lambda,x)=1+\sum\limits_0^4G_s(\lambda,x)+\sum\limits_1^N\frac{R_n(x)e^{i\lambda\zeta_1x}}{\varkappa_n-\lambda}+\sum\limits_1^M\frac{\widehat{R}_m(x)e^{i\lambda\zeta_4x}}{\widehat{\varkappa}_m-\lambda}\label{eq3.46}
\end{equation}
where $R_n(x)$ and $\widehat{R}_m(x)$ are residues at the points $\varkappa_n$ and $\widehat{\varkappa}_m$, and the functions $G_s(\lambda,x)$ corresponding to Cauchy integrals along the rays $\widehat{l}_{\theta\zeta_s}$ ($0\leq s\leq4$) are
$$G_0(\lambda,x)\stackrel{\rm def}{=}\frac1{2\pi i}\int\limits_0^\infty\left\{s_4(t\theta)e^{i\theta t(\zeta_4-\zeta_0)x}\psi_0^+(\theta t,x)+\left(1-\frac{\omega_1}{\omega_0}\right)\psi_1^+(\theta t,x)\right\}\frac{\theta dt}{\theta t-\lambda};$$
$$G_4(\lambda,x)\stackrel{\rm def}{=}\frac1{2\pi i}\int\limits_0^\infty\left\{\frac{\omega_0}{\omega_3}(s_2(\theta t)+s_3(\theta t))e^{i\theta t(\zeta_1-\zeta_4)x}\psi_1^+(\theta t,x)\right.$$
$$\left.+\left(\frac{\omega_2}{\omega_3}e^{i\theta t(\zeta_1-\zeta_3)x}-1\right)\psi_2^+(\theta t,x)\right\}\frac{\theta\zeta_1dt}{\theta\zeta_1t-\lambda};$$
$$G_3(\lambda,x)\stackrel{\rm def}{=}\frac1{2\pi i}\int\limits_0^\infty\left\{\frac{\omega_0}{\omega_3}(s_2(\theta t\zeta_3)+s_3(\theta t\zeta_3))e^{i\theta t(\zeta_0-\zeta_3)x}\psi_3^+(\theta t,x)\right.$$
$$\left.+\left(\frac{\omega_2}{\omega_3}e^{i\theta t(\zeta_1-\zeta_0)x}-1\right)\psi_1^+(\theta t,x)\right\}\frac{\theta\zeta_3dt}{\theta\zeta_3t-\lambda};$$
\begin{equation}
\begin{array}{lll}
{\displaystyle G_2(\lambda,x)\stackrel{\rm def}{=}\frac1{2\pi i}\int\limits_0^\infty\left\{\frac{\omega_0}{\omega_3}[s_2(\theta\zeta_2t)+s_3(\theta\zeta_2t)]e^{it\theta(\zeta_3-\zeta_2)x}\psi_2^+(\theta t,x)\right.}\\
{\displaystyle\left.+\left(\frac{\omega_2}{\omega_3}e^{i\theta t(\zeta_3-\zeta_4)x}-1\right)\psi_2^+(\theta t,x)\right\}\frac{\theta\zeta_2dt}{\theta\zeta_2-\lambda};}\\
{\displaystyle G_1(\lambda,x)\stackrel{\rm def}{=}\frac1{2\pi i}\int\limits_0^\infty\left\{s_1(\theta\zeta_1t)e^{i\theta t(\zeta_2-\zeta_1)x}\psi_0^+(\theta,x)+\left(1-\frac{\omega_4}{\omega_0}\right)\psi_4^+(\theta t,x)\right\}\frac{\theta\zeta_1dt}{\theta\zeta_1t-\lambda}.}
\end{array}\label{eq3.47}
\end{equation}
Inside each sector $\Omega_0$, $\Omega_4^+$, $\Omega_2^+$, $\Omega_3^+$, $\Omega_1^+$, the function $\mathcal{F}(\psi,\lambda,x)$ \eqref{eq3.46} coincides with the functions $\psi_0^+(\lambda,x)$, $\psi_1^+(\lambda,x)$, $\psi_2^+(\lambda,x)$, $\psi_3^+(\lambda,x)$, $\psi_4^+(\lambda,x)$ correspondingly. Therefore, upon calculating boundary values $\lambda\rightarrow\tau\theta\in\widehat{l}_\theta$ ($\lambda\in\Omega_0^-$), $\lambda\rightarrow\tau\theta\zeta_4\in\widehat{l}_{\theta\zeta_4}$ ($\lambda\in\Omega_3^-$), $\lambda\rightarrow\tau\theta\zeta_3\in\widehat{l}_{\theta\zeta_3}$ ($\lambda\in\Omega_2^-$), and, finally, $\lambda\rightarrow\tau\theta\zeta_2\in\widehat{l}_{\theta\zeta_2}$ ($\lambda\in\Omega_1^-$), $\lambda\rightarrow\tau\theta\zeta_1\in\widehat{l}_{\theta\zeta_1}$ ($\lambda\in\Omega_1^-$), and using expressions for the boundary value of the Cauchy-type integral \cite{18,17}, we obtain five equalities:
$${\rm(i)}\,f_1(\tau\theta,x)/\omega_0c=\frac1{\omega_0c}\{\omega_1\psi_1^+(\theta\tau,x)-\omega_0s_4(\theta\tau)e^{it\theta(\zeta_4-\zeta_0)}\psi_0^+(\theta\tau,x)\}$$
$$=1+\sum\limits_0^4G_s(\tau\theta,x)+\sum\limits_1^N\frac{R_n(x)e^{i\theta\tau\zeta_1x}}{\varkappa_n-\tau\theta}+\sum\limits_1^M\frac{\widehat{R}_m(x)e^{i\theta\tau\zeta_4x}}{\widehat{\varkappa_m}-
\tau\theta};$$
$${\rm(ii)}\,\psi_2^+(\tau\theta\zeta_4,x)=1+\sum\limits_0^4G_s(\tau\theta\zeta_4,x)+\sum\limits_1^N\frac{R_n(x)e^{i\theta\tau x}}{\varkappa_n-\tau\theta\zeta_4}+\sum\limits_1^M\frac{\widehat{R}_m(x)e^{i\tau\zeta_3x}}{\widehat{\varkappa}_m-\tau\theta\zeta_4};$$
\begin{equation}
{\rm(iii)}\,\psi_3^+(\tau\theta\zeta_3,x)=1+\sum\limits_0^4G_s(\tau\theta\zeta_3,x)+\sum\limits_1^N\frac{R_n(x)e^{i\theta\tau\zeta_4x}}{\varkappa_n-\tau\theta\zeta_3}+\sum\limits_1^M\frac{\widehat{R}_m(x)e^{i\tau\theta\zeta_2x}}
{\widehat{\varkappa}_m-\tau\theta\zeta_3};\label{eq3.48}
\end{equation}
$${\rm(iv)}\,\psi_4(\tau\theta\zeta_2,x)=1+\sum\limits_0^4G_s(\tau\theta\zeta_2,x)+\sum\limits_1^N\frac{R_n(x)e^{i\tau\zeta_3x}}{\varkappa_n-\tau\theta\zeta_2}+\sum\limits_1^M\frac{\widehat{R}_m(x)}{\widehat{\varkappa}_m-\tau\theta
\zeta_2}e^{i\tau\theta\zeta_1x};$$
$${\rm(v)}\,\psi_1(\tau\theta\zeta_1,x)=1+\sum\limits_0^4G_s(\tau\theta\zeta_1,x)+\sum\limits_1^N\frac{R_n(x)e^{i\tau\theta\zeta_2x}}{\varkappa_n-\tau\theta\zeta_1}+\sum\limits_1^M\frac{\widehat{R}_m(x)}{\widehat{\varkappa}_m-\tau
\theta\zeta_1}e^{i\tau\theta x}.$$
Equalities (i) -- (v) \eqref{eq3.48} give us the linear system of singular integral equations for $\{\psi_k^+(\lambda,x)\}_0^4$. To find $\{R_n(x)\}_1^N$, $\{\widehat{R}_m(x)\}_1^M$, we ought to obtain $N+M$ more equations.

Upon multiplying equality \eqref{eq3.46} by $(\varkappa_p-\lambda)^{-1}$ and integrating over a contour containing $\varkappa_p$, we obtain $N$ equations
\begin{equation}
{\rm(vi)}\,0=1+\sum\limits_0^4G_s(\varkappa_p,x)+\sum\limits_{n\not=p}\frac{R_n(x)e^{i\varkappa_p\zeta_1x}}{\varkappa_n-\varkappa_p}+\sum\limits_1^M\frac{\widehat{R}_m(x)e^{i\varkappa_p\zeta_4x}}{\varkappa_m-\varkappa_p}
\label{eq3.49}
\end{equation}
($0\leq p\leq N$). Besides, we used the fact that coefficients of $(\varkappa_p-\lambda)^0$ vanish because ${\mathcal F}(1,\lambda,x)-1\rightarrow0$ as $\lambda\rightarrow\infty$ ($\lambda\in\Omega_0^-\cap\Omega_2$). Analogously, upon multiplying by $(\widehat{\varkappa}_q-\lambda)^{-1}$ the function ${\mathcal F}(\lambda,x)$ \eqref{eq3.46} and integrating over the contour containing $\widehat{\varkappa}_q$, we obtain $M$ more equations:
\begin{equation}
{\rm(vii)}\,0=1+\sum\limits_0^4G_s(\widehat{\varkappa}_q,x)+\sum\limits_{n=1}^N\frac{R_n(x)e^{i\widehat{\varkappa}_q\zeta_1x}}{\varkappa_n-\widehat{\varkappa}_q}+\sum\limits_{m\not=q}\frac{\widehat{R}_m(x)e^{i\widehat{\varkappa}_q
\zeta_4x}}{\widehat{\varkappa}_m-\varkappa_q}\label{eq3.50}
\end{equation}
($1\leq q\leq M$).

\begin{conclusion}
Relations (i) -- (v) \eqref{eq3.48} and (vi) \eqref{eq3.49}, (vii) \eqref{eq3.50} form the closed linear system of singular integral equations relative to $\{\psi_k^+(\lambda,x)\}_0^4$ and $\{R_n(x)\}_1^N$, $\{R_m(x)\}_1^M$. Scattering coefficients $\{s_k(\lambda)\}_1^4$ and sets of points $\Lambda_1^-$ \eqref{eq3.35}, $\Lambda_4^+$ \eqref{eq3.36} are independent parameters of this system.
\end{conclusion}

\begin{conclusion}
Knowing solution to this system, from formula \eqref{eq3.46} we define function ${\mathcal F}(\psi,\lambda,x)=\psi_1^+(\lambda,x)$ ($\forall\lambda\in\Omega_4^-$). Hence, we know function $\psi_1(\lambda,x)$, then, using Theorem \ref{t2.3} we unambiguously restore the potential $q(x)$ on the right half-axis for all $x\in\mathbb{R}_+$ from the formulas (i) \eqref{eq2.30}, (ii) \eqref{eq2.30}.
\end{conclusion}

\begin{remark}\label{r3.8}
We don't study the problem of solvability of the system (i) \eqref{eq3.48} -- (vii) \eqref{eq3.50} as well as the problem of loosening the condition $v_0(\lambda,x)\in SYMM(\Omega_0)$. For $s_1(\lambda)=s_2(\lambda)=s_3(\lambda)=s_4(\lambda)$, system (i) -- (vii) is simplified and leads to the so-called reflectionless potential which is playing an important role in the soliton theory. The main system of equations (i) \eqref{eq3.48} -- (vii) \eqref{eq3.50} is obtained under the assumption that $v_0(\lambda,x)\in SYMM(\Omega_0)$. The problem of loosening or abandonment this condition remains unsolved.
\end{remark}

\section{Scattering problem for the waves incident from ``$-\infty$''}\label{s4}

{\bf 4.1} Analogously to \eqref{eq3.1}, we expand every Jost solution $\{v_k(\lambda,x)\}_0^4$ by the other Jost solutions $\{u_k(\lambda,x)\}_0^4$,
\begin{equation}
v_k(\lambda,x)=\sum\limits_{s=0}^4\widetilde{t}_{k,s}(\lambda)u_s(\lambda,x)\quad(0\leq k\leq4)\label{eq4.1}
\end{equation}
or, in the matrix form,
\begin{equation}
v(\lambda,x)=\widetilde{T}(\lambda)u(\lambda,x)\label{eq4.2}
\end{equation}
where $v(\lambda,x)\stackrel{\rm def}{=}\col[v_0(\lambda,x),...,v_4(\lambda,x)]$; $u(\lambda,x)\stackrel{\rm def}{=}\col[u_0(\lambda,x),...,u_4(\lambda,x)]$, and $\widetilde{T}(\lambda)$ is the dual transition matrix:
\begin{equation}
\widetilde{T}(\lambda)\stackrel{\rm def}{=}\left[
\begin{array}{cccc}
\widetilde{t}_{0,0}(\lambda),&...,&\widetilde{t}_{0,4}(\lambda)\\
\widetilde{t}_{1,0}(\lambda),&...,&\widetilde{t}_{1,4}(\lambda)\\
\cdots&\cdots&\cdots\\
\widetilde{t}_{4,0}(\lambda),&...,&\widetilde{t}_{4,4}(\lambda)
\end{array}\right].\label{eq4.3}
\end{equation}

\begin{remark}\label{r4.1}
For $k=0$, equalities \eqref{eq4.1} yield
\begin{equation}
v_0(\lambda,x)=\sum\limits_{s=0}^4\widetilde{t}_{0,s}(\lambda)u_s(\lambda,x).\label{eq4.4}
\end{equation}
The other equalities in \eqref{eq4.1} follow from \eqref{eq4.4}, in view of \eqref{eq2.18} upon the substitution $\lambda\rightarrow\lambda\zeta_k$, besides, for $\widetilde{t}_{k,s}(\lambda)$, equalities analogous to \eqref{eq3.5} hold.

Equations \eqref{eq3.2} and \eqref{eq4.2} imply that $T(x)\widetilde{T}(\lambda)=I$. Taking into account that $T(\lambda)\cdot JT^*(\overline{\lambda})J^*\overline{\eta}$, due to \eqref{eq3.21} (see Theorem \ref{t3.2}), we have
\begin{equation}
\widetilde{T}(\lambda)=\overline{\eta}JT^*(\overline{\lambda})J,\label{eq4.5}
\end{equation}
besides, $J$ is given by \eqref{eq3.20}.
\end{remark}

Equation \eqref{eq4.4} implies the equality
\begin{equation}
\widetilde{r}_0(\lambda)v_0(\lambda,x)=u_0(\lambda,x)+\sum\limits_1^4\widetilde{s}_l(\lambda)u_l(\lambda,x)\label{eq4.6}
\end{equation}
where
\begin{equation}
\widetilde{r}_0(\lambda)\stackrel{\rm def}{=}1/\widetilde{t}_{0,0}(\lambda);\quad\widetilde{s}_l(\lambda)\stackrel{\rm def}{=}\widetilde{t}_{0l}(\lambda)/\widetilde{t}_{0,0}(\lambda)\quad(1\leq l\leq4)\label{eq4.7}
\end{equation}

\begin{remark}\label{r4.2}
Equation (2.4.b) imply the asymptotic at ``$-\infty$'',
$$\widetilde{r}_0(\lambda)u_0(\lambda,x)=e^{i\lambda\zeta_0x}+\sum\limits_1^4\widetilde{s}_l(\lambda)e^{i\lambda\zeta_lx}\quad(x\rightarrow-\infty)$$
and thus the incident from ``$-\infty$'' wave $e^{i\lambda\zeta_0}$ has the {\bf reflected (scattered) wave} ${\displaystyle\sum\limits_1^4\widetilde{s}_l(\lambda)e^{i\lambda\zeta_lx}}$ where $\widetilde{s}_l(\lambda)$ are {\bf scattering coefficients}. At ``$+\infty$'', function $r_0(\lambda)v_0(\lambda,x)$ has the following asymptotic behavior, in view of (2.4.a),
$$\widetilde{r}_0(\lambda)v_0(\lambda,x)\rightarrow r_0(\lambda)e^{i\lambda\zeta_0x}\quad(x\rightarrow\infty)$$
It is natural to consider $\widetilde{r}_0(\lambda)$ to be the {\bf transition coefficient} of the wave $e^{i\lambda\zeta_0x}$ incident from ``$-\infty$''.
\end{remark}

Upon substituting $\lambda\rightarrow\lambda\zeta_k$ and using \eqref{eq2.18}, \eqref{eq4.6} implies the following scattering problems for the incident from ``$-\infty$'' waves $e^{i\lambda\zeta_kx}$:
\begin{equation}
\widetilde{r}_k(\lambda)v_k(\lambda,x)=u_k(\lambda,x)+\sum\limits_{l\not=k}\widetilde{s}_{k,l}(\lambda)u_l(\lambda,x)\quad(1\leq k\leq4)\label{eq4.8}
\end{equation}
where
\begin{equation}
\widetilde{r}_k(\lambda)\stackrel{\rm def}{=}1/\widetilde{t}_{k,k}(\lambda);\quad\widetilde{s}_{k,l}(\lambda)=\widetilde{t}_{k,l}(\lambda)/\widetilde{t}_{k,k}(\lambda)\quad(l\not=k).\label{eq4.9}
\end{equation}
Note that transition coefficients $\{\widetilde{r}_k(\lambda)\}_0^4$, as well as the scattering coefficients $\{\widetilde{s}_{k,l}(\lambda)\}_0^4$ \eqref{eq4.7}, \eqref{eq4.9}, are expressed via transition coefficients $\{t_k(\lambda)\}_0^4$ and scattering coefficients $\{s_{k,l}(\lambda)\}_0^4$ of the scattering problem for the waves incident from ``$+\infty$'' in view of \eqref{eq4.5}.

Analogously to \eqref{eq3.37}, write equalities \eqref{eq4.8} as
$${\rm(i)}\,\widetilde{r}_1(\lambda)v_1(\lambda,x)=\widetilde{s}_4(\lambda\zeta_1)u_0(\lambda,x)+u_1(\lambda,x)+\widetilde{s}_1(\lambda\zeta_1)u_2(\lambda,x)+\widetilde{s}_2(\lambda\zeta_1)u_3(\lambda,x)$$
$$+\widetilde{s}_3(\lambda\zeta_1)u_4(\lambda,x);$$
$${\rm(ii)}\,\widetilde{r}_2(\lambda)v_2(\lambda,x)=\widetilde{s}_3(\lambda\zeta_2)u_0(\lambda,x)+\widetilde{s}_4(\lambda\zeta_2)u_1(\lambda,x)+u_2(\lambda,x)+\widetilde{s}_1(\lambda\zeta_2)u_3(\lambda,x)$$
\begin{equation}
+\widetilde{s}_2(\lambda\zeta_2)u_4(\lambda,x);\label{eq4.10}
\end{equation}
$${\rm(iii)}\,\widetilde{r}_3(\lambda)v_3(\lambda,x)=\widetilde{s}_2(\lambda\zeta_3)u_0(\lambda,x)+\widetilde{s}_3(\lambda\zeta_3)u_1(\lambda,x)+\widetilde{s}_4(\lambda\zeta_3)u_2(\lambda,x)$$
$$+u_3(\lambda,x)+\widetilde{s}_2(\lambda\zeta_2)u_4(\lambda,x);$$
$${\rm(iv)}\,\widetilde{r}_4(\lambda)v_4(\lambda,x)=\widetilde{s}_1(\lambda\zeta_4)u_0(\lambda,x)+\widetilde{s}_2(\lambda\zeta_4)u_1(\lambda,x)+s_3(\lambda\zeta_4)u_2(\lambda,x)$$
$$+\widetilde{s}_4(\lambda\zeta_4)u_3(\lambda,x)+u_4(\lambda,x),$$
besides, $\widetilde{r}_k(\lambda)\stackrel{\rm def}{=}\widetilde{r}_0(\lambda\zeta_k)$ ($0\leq k\leq4$).
\vspace{5mm}

{\bf 4.2} Calculating Wronskians $w(v_0,u_1,u_2,u_3)$ and $w(v_0,v_2,v_3,v_4)$, where $v_0(\lambda,x)$ is given by \eqref{eq4.6}, we have
$${\rm(i)}\,\widetilde{r}_0(\lambda)w(v_0,u_1,u_2,u_3)=w(u_0,u_1,u_2,u_3)-\widetilde{s}_4(\lambda)w(u_1,u_2,u_3,u_4);$$
$${\rm(ii)}\,\widetilde{r}_0(\lambda)w(v_0,u_2,u_3,u_4)=w(u_0,u_2,u_3,u_4)+\widetilde{s}_1(\lambda)w(u_1,u_2,u_3,u_4).$$
Now using analogues of the formulas \eqref{eq3.12} for the Wronskians $W_{k,s}(u,\lambda,x)$, we obtain
$${\rm(i)}\,\widetilde{r}_0(\lambda)w(v_0,u_1,u_2,u_3)=(i\lambda)^6\{\omega_1u_1^+(\lambda,x)-\widetilde{s}_4(\lambda)\omega_0u_0^+(\lambda,x)\};$$
$${\rm(ii)}\,\widetilde{r}_0(\lambda)w(v_0,u_2,u_3,u_4)=(i\lambda)^6\{\omega_4u_4^+(\lambda,x)+\widetilde{s}_1(\lambda)\omega_0u_0^+(\lambda,x)\}.$$
Assuming that $u_k(\lambda,x)=\varphi_k(\lambda,x)e^{i\lambda\zeta_kx}$, hence we find
\begin{equation}
\begin{array}{lll}
{\rm(i)}\,g_1(\lambda,x)-\omega_1\varphi_1^+(\lambda,x)=-\omega_0s_4(\lambda)e^{i\lambda(\zeta_4-\zeta_0)x}\varphi_0^+(\lambda,x)&(l_{\theta\zeta_1});\\
{\rm(ii)}\,g_4(\lambda,x)-\omega_4\varphi_4^+(\lambda,x)=\omega_0s_0(\lambda)e^{i\lambda(\zeta_1-\zeta_0)x}\varphi_0^+(\lambda,x)&(l_\theta)
\end{array}\label{eq4.11}
\end{equation}
where
\begin{equation}
\begin{array}{cccc}
g_1(\lambda,x)\stackrel{\rm def}{=}\widetilde{r}_0(\lambda)w(v_0,u_1,u_2,u_3)e^{i\lambda\zeta_4x}(i\lambda)^{-6};\\
g_4(\lambda,x)\stackrel{\rm def}{=}\widetilde{r}_0(\lambda)w(v_0,u_2,u_3,u_4)e^{i\lambda\zeta_1x}(i\lambda)^{-6}.
\end{array}\label{eq4.12}
\end{equation}
Repeating considerations 3.8, it is easy to make sure that functions $g_1(\lambda,x)$, $g_4(\lambda,x)$ \eqref{eq4.12} are holomorphic in the sectors $\Omega_0\cap\Omega_3^-$ and $\Omega_0\cap\Omega_2^-$ correspondingly, besides,
$$g_1(\lambda,x)\rightarrow c\omega_0\quad(\lambda\rightarrow\infty,\lambda\in\Omega_0\cap\Omega_3^-);\quad g_4(\lambda,x)\rightarrow c\omega_0\quad(\lambda\rightarrow\infty,\lambda\in\Omega_0\cap\Omega_2^-).$$
As a result, we have two Riemann boundary value problems: (4.12.i) on the ray $l_{\theta\zeta_1}$ for holomorphic in the sectors $\Omega_1$ and $\Omega_0\cap\Omega_3^-$ functions $\varphi_1^+(\lambda,x)$ and $g_1(\lambda,x)$; (4.12.ii) on the ray $l_\theta$ for holomorphic in the sectors $\Omega_4$ and $\Omega_0\cap\Omega_2^-$ functions $\varphi_4(\lambda,x)$ and $g_4(\lambda,x)$.

\begin{picture}(200,200)
\put(0,100){\line(1,0){200}}
\put(100,0){\vector(0,1){200}}
\put(0,133){\vector(3,-1){200}}
\put(200,133){\vector(-3,-1){200}}
\put(105,150){$g_4(\lambda,x)$}
\qbezier(100,140)(135,145)(141,115)
\qbezier(140,85)(155,105)(143,116)
\qbezier(56,86)(50,100)(58,115)
\qbezier(100,141)(65,150)(57,116)
\put(0,105){$\varphi_1^+(\lambda,x)$}
\put(55,102){$\Omega_1$}
\put(20,145){$g_1(\lambda,x)$}
\put(170,110){$\varphi_4^+(\lambda,x)$}
\put(130,102){$\Omega_4$}
\put(105,40){$u_0^+(\lambda,x)$}
\end{picture}

\hspace{20mm} Fig. 8. Boundary value problems \eqref{eq4.11}
\vspace{5mm}

{\bf 4.3} Using representation \eqref{eq4.6} of function $v_0(\lambda,x)$, calculate Wronskians $w(v_0,u_1,u_2,u_4)$ and $w(v_0,u_1,u_3,u_4)$ and using Lemma \ref{l3.1}, we have
$${\rm(i)}\,\widetilde{r}_0(\lambda)w(v_0,u_1,u_2,u_4)=w(u_0,u_1,u_2,u_4)+\widetilde{s}_3(\lambda)w(u_1,u_2,u_3,u_4);$$
$${\rm(ii)}\,\widetilde{r}_0(\lambda)w(v_0,u_1,u_3,u_4)=w(u_0,u_1,u_3,u_4)-\widetilde{s}_2(\lambda)w(u_1,u_2,u_3,u_4)$$
and using analogues of the formulas \eqref{eq3.12} for $w_{k,s}(u,\lambda,x)$, we obtain
$${\rm(i)}\,\widetilde{r}_0(\lambda)w(v_0,u_1,u_2,u_4)=(i\lambda)^6\{\omega_2u_2^+(\lambda,x)+\widetilde{s}_3(\lambda)\omega_0u_0^+(\lambda,x)\};$$
$${\rm(ii)}\,\widetilde{r}_0(\lambda)w(v_0,u_1,u_3,u_4)=(i\lambda)^6\{\omega_3u_3^+(\lambda,x)-\widetilde{s}_2(\lambda)\omega_0u_0^+(\lambda,x)\}.$$
Assuming that $u_0(\lambda,x)\in SYMMM(\Omega_0^-)$, i.e., $u_0(\theta t,x)=u_0(\theta\zeta_1,x)$ ($\forall t\in\mathbb{R}_+$), we have $\left.u_2(\lambda,x)\right|_{l_{\theta\zeta_3}}=\left.u_3(\lambda,x)\right|_{l_{\theta\zeta_3}}$. Therefore, upon subtracting two last equalities, we have
\begin{equation}
\omega_3u_3^+(\lambda,x)-\omega_2u_2^+(\lambda,x)=[\widetilde{s}_2(\lambda)+\widetilde{s}_3(\lambda)]\omega_0u_0^+(\lambda,x)\quad(\forall\lambda\in l_{\theta\zeta_3}).\label{eq4.13}
\end{equation}
As a result, we arrive at the Riemann problem on ray $l_{\theta\zeta_3}$ for the functions $u_2^+(\lambda,x)$ and $u_3^+(\lambda,x)$ holomorphic in the sectors $\Omega_2$ and $\Omega_3$ (see Fig.)

\begin{picture}(200,200)
\put(0,100){\line(1,0){200}}
\put(100,0){\vector(0,1){200}}
\put(0,133){\vector(3,-1){200}}
\put(200,133){\vector(-3,-1){200}}
\qbezier(100,60)(135,70)(144,85)
\qbezier(100,58)(70,60)(57,82)
\put(105,20){$l_{\theta\zeta_3}$}
\put(35,40){$u_2^+(\lambda,x)$}
\put(105,40){$u_3^+(\lambda,x)$}
\put(75,70){$\Omega_2$}
\put(110,72){$\Omega_3$}
\end{picture}

\hspace{20mm} Fig. 9. Boundary value problem \eqref{eq4.13}

Upon substitutions $\lambda\rightarrow\lambda\zeta_1$, $\lambda\rightarrow\lambda\zeta_4$ in \eqref{eq4.13}, we obtain the following boundary value problems on $l_{\theta\zeta_4}$ and $l_{\theta\zeta_2}$:
\begin{equation}
\begin{array}{ccc}
{\rm(i)}\,\omega_3u_2^+(\lambda,x)-\omega_2u_1^+(\lambda,x)=[\widetilde{s}_2(\lambda\zeta_1)+\widetilde{s}_3(\lambda\zeta_1)]\omega_0u_4^+(\lambda,x)\quad(\forall\lambda\in l_{\theta\zeta_4});\\
{\rm(ii)}\,\omega_3u_4^+(\lambda,x)-\omega_2u_3^+(\lambda,x)=[\widetilde{s}_2(\lambda\zeta_4)+\widetilde{s}_3(\lambda\zeta_4)]\omega_0u_1^+(\lambda,x)\quad(\forall\lambda\in l_{\theta\zeta_2}).
\end{array}\label{eq4.14}
\end{equation}
\vspace{5mm}

{\bf 4.4} Define the following holomorphic function:
\begin{equation}
\widetilde{\mathcal F}(\varphi,\lambda,x)\stackrel{\rm def}{=}\left\{
\begin{array}{llllll}
g_4(\lambda,x)/c\omega_0;\quad(\lambda\in\Omega_0\cap\Omega_3^-);\\
g_1(\lambda,x)/c\omega_0;\quad(\lambda\in\Omega_0\cap\Omega_0^-);\\
\varphi_1^+(\lambda,x);\quad\lambda\in\Omega_1;\\
\varphi_2^+(\lambda,x);\quad\lambda\in\Omega_2;\\
\varphi_3^+(\lambda,x);\quad\lambda\in\Omega_3;\\
\varphi_4^+(\lambda,x);\quad\lambda\in\Omega_4.
\end{array}\right.\label{eq4.15}
\end{equation}
Function $\widetilde{\mathcal F}(\varphi,\lambda,x)$ is holomorphic in the sectors $\{\Omega_k\}_1^4$ and meromorphic in the sectors $\Omega_0\cap\Omega_3^-$, $\Omega_0\cap\Omega_2^-$, besides, $\widetilde{\mathcal F}(\varphi,\lambda,x)\rightarrow1$ ($\lambda\rightarrow\infty$) inside each of the sectors. Note that in the sectors $\Omega_0\cap\Omega_3^-$ and $\Omega_0\cap\Omega_2^-$ function $\widetilde{\mathcal{F}}(\varphi,\lambda,x)$ has simple poles which coincide with zeros of function $t_{0,0}^+(\lambda)$.

Equations (4.11.i), (4.11.ii) and \eqref{eq4.13} and \eqref{eq4.14} yield the following jump problems:
$${\rm(a)}\,\varphi_1^+(\lambda,x)-g_1(\lambda,x)=\widetilde{s}_4(\lambda)e^{i\lambda(\zeta_4-\zeta_0)x}\varphi_0^+(\lambda,x)+\left(1-\frac{\omega_1}{\omega_0}\right)\varphi_1^+(\lambda,x)\quad(l_\theta);$$
$${\rm(b)}\,\varphi_4^+(\lambda,x)-g_4(\lambda,x)=-\widetilde{s}_1(\lambda)e^{i\lambda(\zeta_1-\zeta_0)x}\varphi_0^+(\lambda,x)+\left(1-\frac{\omega_4}{\omega_0}\right)\varphi_4^+(\lambda,x)\quad(l_{\theta\zeta_1});$$
\begin{equation}
{\rm(c)}\,\varphi_3^+(\lambda,x)-\varphi_2^+(\lambda,x)=\frac{\omega_0}{\omega_3}[\widetilde{s}_2(\lambda)+\widetilde{s}_3(\lambda)]e^{i\lambda(\zeta_2-\zeta_0)x}\varphi_0^+(\lambda,x)\label{eq4.16}
\end{equation}
$$+\left(\frac{\omega_2}{\omega_3}e^{i\lambda(\zeta_2-\zeta_3)x}-1\right)\varphi_2^+(\lambda,x)\quad(l_{\theta\zeta_3});$$
$${\rm(d)}\,\varphi_2^+(\lambda,x)-\varphi_1^+(\lambda,x)=\frac{\omega_0}{\omega_3}[\widetilde{s}_2(\lambda\zeta_1)+\widetilde{s}_3(\lambda\zeta_1)]e^{i\lambda(\zeta_3-\zeta_1)x}\varphi_4^+(\lambda,x)$$
$$+\left(\frac{\omega_2}{\omega_3}e^{i\lambda(\zeta_3-\zeta_4)x}-1\right)\varphi_1^+(\lambda,x)\quad(l_{\theta\zeta_4});$$
$${\rm(e)}\,\varphi_4^+(\lambda,x)-\varphi_3^+(\lambda,x)=\frac{\omega_0}{\omega_3}[\widetilde{s}_2(\lambda\zeta_4)+\widetilde{s}_3(\lambda\zeta_4)]e^{i\lambda(\zeta_1-\zeta_4)x}\varphi_1^+(\lambda,x)$$
$$+\left(\frac{\omega_2}{\omega_3}e^{i\lambda(\zeta_1-\zeta_2)x}-1\right)\varphi_3^+(\lambda,x)\quad(l_{\theta\zeta_2})$$
Calculate function $\widetilde{\mathcal{F}}(\varphi,\lambda,x)$ \eqref{eq4.14} using the Cauchy integral
$$\widetilde{\mathcal{F}}(\varphi,\lambda,x)=1+\frac1{2\pi i}\cdot\int\limits_{\widetilde{\gamma}}\frac{\mathcal{F}(\varphi,w,x)-1}{w-\lambda}dw$$
where $\widetilde{\gamma}$ is the contour obtained from contour $\gamma$ (see Fig. 7). Then, taking into account \eqref{eq4.15}, we obtain the following expression for function $\widetilde{\mathcal{F}}(\varphi,\lambda,x)$,
\begin{equation}
\widetilde{F}(\varphi,\lambda,x)=1+\sum\limits_{s=0}^4\widetilde{G}_s(\lambda,x)+\sum\limits_1^N\frac{\widetilde{R}_n(x)e^{i\lambda\zeta_4x}}{\overline{\varkappa}_n-\lambda}+\sum\limits_1^M\frac{\widetilde{\widehat{R}}_m(\lambda)
e^{i\lambda\zeta_1x}}{\overline{\widehat{\varkappa}}_m-\lambda}\label{eq4.17}
\end{equation}
where $\widetilde{R}_n(\lambda,x)$ and $\widetilde{\widehat{R}}_m(x)$ are residues in the points $\overline{\varkappa}_n$ and $\overline{\widehat{\varkappa}}_m$ of the functions $g_1(\lambda,x)$ and $g_4(\lambda,x)$ on the rays $l_{\theta\zeta_1}$ and $l_{\theta\zeta_4}$. For the functions $\widetilde{G}_s(\lambda,x)$ corresponding to the Cauchy integrals along the rays $l_{\theta\zeta_s}$ ($0\leq s\leq4$), we have
$$\widetilde{G}_0(\lambda,x)\stackrel{\rm def}{=}\frac1{2\pi i}\int\limits_0^\infty\left\{\widetilde{s}_4(t\theta)e^{i\theta t(\zeta_4-\zeta_0)x}\varphi_0^+(\theta t,x)+\left(1-\frac{\omega_1}{\omega_0}\right)\varphi_1^+(\theta t,x)\right\}\frac{\theta dt}{\theta t-\lambda};$$
$$\widetilde{G}_4(\lambda,x)\stackrel{\rm def}{=}\frac1{2\pi i}\int\limits_0^\infty\left\{\frac{\omega_0}{\omega_3}[\widetilde{s}_2(\theta t)+s_3(\theta t)]e^{i\theta t(\zeta_1-\zeta_4)x}\varphi_1^+(\theta t,x)\right.$$
$$\left.+\left(\frac{\omega_3}{\omega_0}e^{i\theta t(\zeta_1-\zeta_3)x}-1\right)\varphi_2^+(\lambda,x)\right\}\frac{\theta\zeta_1dt}{\theta\zeta_1t-\lambda};$$
$$\widetilde{G}_3(\lambda,x)\stackrel{\rm def}{=}\frac1{2\pi i}\int\limits_0^\infty\left\{\frac{\omega_0}{\omega_3}[\widetilde{s}_2(\theta t\zeta_3)+\widetilde{s}_3(\theta t\zeta_3)]e^{i\theta t(\zeta_0-\zeta_3)x}\varphi_3^+(\theta t,x)\right.$$
\begin{equation}
\left.+\left(\frac{\omega_2}{\omega_3}e^{i\theta t(\zeta_1-\zeta_0)x}-1\right)\varphi_1^+(\theta t,x)\right\}\frac{\theta\zeta_3dt}{\theta\zeta_3t-\lambda};\label{eq4.18}
\end{equation}
$$\widetilde{G}_2(\lambda,x)\stackrel{\rm def}{=}\frac1{2\pi i}\int\limits_0^\infty\left\{\frac{\omega_0}{\omega_3}[s_2(\theta\zeta_2t)+s_3(\theta\zeta_2t)]e^{i\theta t(\zeta_3-\zeta_2)x}\varphi_2^+(\lambda,x)\right.$$
$$\left.+\left(\frac{\omega_2}{\omega_3}e^{i\theta t(\zeta_3-\zeta_4)x}-1\right)\varphi_2^+(\lambda,x)\right\}\frac{\theta\zeta_2dt}{\theta\zeta_2t-\lambda};$$
$$\widetilde{G}_1(\lambda,x)\stackrel{\rm def}{=}\frac1{2\pi i}\int\limits_0^\infty\left\{\widetilde{s}_1(\theta\zeta,t)e^{i\theta t(\zeta_2-\zeta_1)x}\varphi_0^+(\lambda,x)+\left(1-\frac{\omega_4}{\omega_0}\right)\varphi_4^+(\lambda,x)\right\}\frac{\theta\zeta_1dt}{\theta\zeta_1t-\lambda}$$
inside the sectors $\Omega_4$, $\Omega_2$, $\Omega_3$, $\Omega_1$ for function $\widetilde{\mathcal{F}}(\varphi,\lambda,x)$ \eqref{eq4.16}, which corresponds to the functions $\varphi_4^+(\lambda,x)$, $\varphi_2^+(\lambda,x)$, $\varphi_3^+(\lambda,x)$, $\varphi_1^+(\lambda,x)$ respectively. Calculating the boundary values as $\lambda\rightarrow\tau\theta\in l_\theta$ ($\lambda\in\Omega_0$); $\lambda\rightarrow\tau\theta\zeta_4\in l_{\theta\zeta_4}$ ($\lambda\in\Omega_3$); $\lambda\rightarrow\tau\theta\zeta_3\in l_{\theta\zeta_3}$ ($\lambda\in\Omega_2$), and, finally, $\lambda\rightarrow\tau\theta\zeta_2\in l_{\theta\zeta_2}$ ($\lambda\in\Omega_1$); $\lambda\rightarrow\tau\theta\zeta_1$ ($\lambda\in\Omega_1$) and using boundary values of the Cauchy type integral \cite{17,18}, we obtain the five equations
$${\rm(i)}\,g_1(\theta\tau,x)/\omega_0c=\frac1{\omega_0c}\{\omega_1\varphi_1^+(\tau\theta,x)-\omega_0\widetilde{s}_4(\theta\tau)e^{i\theta\tau(\zeta_4-\zeta_0)x}\}$$
$$=1+\sum\limits_{s=0}^4\widetilde{G}_s(\tau\theta,x)+\sum\limits_1^N\frac{\widetilde{R}_n(\lambda)e^{i\theta\tau\zeta_4x}}{\overline{\varkappa}_m-\tau\theta}+\sum\limits_1^M\frac{\widetilde{\widehat{R}}e^{i\theta t\zeta_1x}}{\overline{\widehat{\varkappa}}_m-\theta\tau};$$
$${\rm(ii)}\,\varphi_4^+(\tau\theta\zeta_4,x)=1+\sum\limits_0^4\widetilde{G}_s(\tau\theta\zeta_4,x)+\sum\limits_1^N\frac{\widetilde{R}_n(x)e^{it\theta\zeta_3x}}{\overline{\varkappa}_n-\tau\theta\zeta_4}+\sum\limits_1^M\frac{
\widetilde{\widehat{R}}_m(x)e^{i\theta tx}}{\overline{\widehat{\varkappa}}-\tau\theta\zeta_4};$$
\begin{equation}
{\rm(iii)}\,\varphi_3^+(\tau\theta\zeta_3,x)=1+\sum\limits_0^4\widetilde{G}_s(\tau\theta\zeta_3,x)+\sum\limits_1^N\frac{\widetilde{R}_n(x)e^{it\zeta_2x}}{\overline{\varkappa}_n-\tau\theta\zeta_3}+\sum\limits_1^M\frac{\widetilde{
\widehat{R}}_m(x)e^{i\theta t\zeta_4x}}{\overline{\widehat{\varkappa}}-\tau\theta\zeta_3};\label{eq4.19}
\end{equation}
$${\rm(iv)}\,\varphi_2^+(\tau\theta\zeta_2,x)=1+\sum\limits_0^4\widetilde{G}_s(\tau\theta\zeta_2,x)+\sum\limits_1^N\frac{\widetilde{R}_n(x)e^{it\zeta_1x}}{\overline{\varkappa}_n-\tau\theta\zeta_2}+\sum\limits_1^M\frac{\widetilde{
\widehat{R}}_m(x)e^{it\theta\zeta_3x}}{\overline{\widehat{\varkappa}}-\tau\theta\zeta_2};$$
$${\rm(v)}\,\varphi_1^+(\tau\theta\zeta_1,x)=1+\sum\limits_0^4\widetilde{G}_s(\tau\theta\zeta_1,x)+\sum\limits_1^N\frac{\widetilde{R}_n(x)e^{it\theta x}}{\overline{\varkappa}_n-\tau\theta\zeta_1}+\sum\limits_1^M\frac{\widetilde{\widehat{R}}_m(x)e^{it\theta\zeta_2x}}{\overline{\widehat{\varkappa}}-\tau\theta\zeta_1}.$$
Equalities (i) -- (v) \eqref{eq4.18} form a system of linear singular equations for the functions $\{\varphi_k^+(\lambda,x)\}_0^4$. In order to find $\{\widetilde{R}_n(x)\}_1^N$ and $\{\widetilde{\widehat{R}}_m(x)\}_1^M$, we obtain the $N+M$ more equations. Upon multiplying equality \eqref{eq4.16} by $(\overline{\varkappa}_p-\lambda)^{-1}$ and integrating along the contour containing point $\varkappa_p$, we obtain the $N$ equations
\begin{equation}
0=1+\sum\limits_0^4\widetilde{G}_s(\overline{\varkappa}_p,x)+\sum\limits_{n\not=p}\frac{\widetilde{R}_n(x)e^{i\overline{\varkappa}_n\zeta_4x}}{\overline{\varkappa}_n-\overline{\varkappa}_p}+\sum\limits_1^M\frac{\widetilde{
\widehat{R}}_m(x)e^{i\overline{\varkappa}_p\zeta_1x}}{\overline{\varkappa}_m-\overline{\varkappa}_q}\label{eq4.20}
\end{equation}
($1\leq p\leq N$).
Analogously, multiplying \eqref{eq4.17} by $(\overline{\widehat{\varkappa}}_q-\lambda)^{-1}$ and integrating along the contour containing $\overline{\widehat{\varkappa}}_q$, we obtain $M$ more equations.
\begin{equation}
0=1+\sum\limits_0^4\widetilde{G}_s(\overline{\widehat{\varkappa}}_q,x)+\sum\limits_1^N\frac{\widetilde{R}_n(x)e^{i\overline{\widehat{\varkappa}}_q\zeta_4x}}{\overline{\varkappa}_n-\overline{\widehat{\varkappa}}_q}+\sum\limits_{
m\not=q}\frac{\widetilde{\widehat{R}}_m(x)e^{i\overline{\widehat{\varkappa}}_q\zeta_1x}}{\overline{\widehat{\varkappa}}_m-\overline{\widehat{\varkappa}}_q}\label{eq4.21}
\end{equation}
($1\leq q\leq M$).

\begin{conclusion}
Equalities (i) -- (v) \eqref{eq4.18}, \eqref{eq4.19}, \eqref{eq4.20} represent a closed linear system of integral singular equations relative to $\{\varphi_k^+(\lambda,x)\}_0^4$ and $\{\widetilde{R}_n(x)\}_1^N$, $\{\widetilde{\widehat{R}}_m(x)\}_1^M$. Scattering coefficients $\{\widetilde{s}_k(\lambda)\}_1^4$ and points $\overline{\Lambda}_1^-$; $\overline{\Lambda}_4^+$ are the independent parameters of this system ($\Lambda_1^-$ and $\Lambda_4^+$ are given by \eqref{eq3.35}, \eqref{eq3.36}).
\end{conclusion}

\begin{conclusion}
Knowing solution to the system \eqref{eq4.18} -- \eqref{eq4.20}, using formula \eqref{eq4.16}, we construct function $\widetilde{\mathcal{F}}(\varphi,\lambda,x)=\varphi_1^+(\lambda,x)$ ($\lambda\in\Omega_1$). From function $\varphi_1^+(\lambda,x)$, we define $\varphi(\lambda,x)$, and afterwards, using Theorem \ref{t2.4}, from formulas \eqref{eq2.33}, we unambiguously find potential $q(x)$ on the left half-axis $x\in\mathbb{R}_-$.
\end{conclusion}

\begin{remark}
Requirement $u_0(\lambda,x)\in SYMM(\Omega_0^-)$ supposes some symmetry, as well as the condition $v_0(\lambda,x)\in SYMM(\Omega_0)$, and is additional. Problem of solvability of the equation systems \eqref{eq3.48} -- \eqref{eq3.50} and \eqref{eq4.18} -- \eqref{eq4.20} is a subject of another study and we don't consider it here.
\end{remark}

{\bf Acknowledgement.} The research was /partially/ supported by NASU project M 22-13 and by the Akhiezer Foundation.

\renewcommand{\refname}{ \rm 
\end{Large}
\end{document}